\documentstyle[12pt]{amsart}

\newtheorem{thm}{Theorem}
\newtheorem{cor}{Corollary}
\newtheorem{lem}{Lemma}
\newtheorem{defn}{Definition}

\theoremstyle{definition}
\newtheorem{rem}{Remark}
\newtheorem{ex}{Example}

\begin{document}

\title[Picard Modular Group]{An explicit fundamental domain for the Picard modular group in two complex dimensions}         
\author{G\'abor Francsics}
\address{Department of Mathematics\\Michigan State University}       
\author{Peter D. Lax}
\address{Courant Institute\\New York University}
\date{July 11, 2005The first author is grateful to the Mathematical Institute of University of Oxford for its hospitality. The research of the first author was partially supported by a Michigan State University IRGP grant. }

\begin{abstract}
Our main goal in this paper is to construct the first explicit fundamental domain of the Picard modular group acting on the complex hyperbolic space  ${\bf CH}^{2}$. The complex hyperbolic space is a Hermitian symmetric space, its bounded realization is the unit ball in ${\bf C}^2$ equipped with the Bergman metric. The Picard modular group is a discontinuous holomorphic automorphism subgroup of $SU(2,1)$ with Gaussian integer entries. This fundamental domain has finite volume, one cusp, explicitly given boundary surfaces and an interesting symmetry.  
\end{abstract}
\maketitle

\section{Introduction}
Our main goal in this paper is to construct the first explicit fundamental domain of the Picard modular group $PU(2,1;{\bf Z}[i])$ acting on the complex hyperbolic space ${\bf CH}^{2}$. The complex hyperbolic space ${\bf CH}^{n}$ is a rank one Hermitian symmetric space of noncompact type. Its bounded realization is the complex unit ball of ${\bf C}^{n}$ equipped with the Bergman metric. The Picard modular group $PU(n,1;{\bf Z}[i])$ is a discontinuous holomorphic automorphism subgroup of ${\bf CH}^{n}$ with Gaussian integer entries. It is the higher dimensional analogue of the modular group, $PSL(2,{\bf Z})$, in ${\bf C}^{n}$. 

Fundamental domains for lattices in rank one symmetric spaces attracted much attention during the last three decades. Although remarkable progress has been achieved, several important problems related to arithmeticity, existence of embedded eigenvalues in the continuous spectrum etc., are still open. The general structure of a fundamental domain for lattices is well known since the work of Garland-Raghunathan \cite{GR}, for example. However there are very few fundamental domains known completely explicitly. This is especially true for complex hyperbolic spaces. Constructing explicit fundamental domains in complex hyperbolic spaces is much more difficult than in real hyperbolic spaces. This phenomenon is well known since the work of Mostow \cite{M}. Recently very strong progress has been made in constructing explicit fundamental domains for discrete subgroups of complex hyperbolic spaces; see for example, the work of Cohen \cite{C}, Holzapfel \cite{H1}, \cite{H2}, Goldman \cite{G}, Goldman-Parker \cite{GP}, Falbel-Parker \cite{FP1}, \cite{FP2}, Schwartz \cite{Sch}, Francsics-Lax \cite{FL1}, \cite{FL2}.  However, explicit fundamental domains do not seem to be known in the literature for the Picard modular groups, e.g., see the comment in \cite{FP2}. 

Our method builds on the construction of a semi-explicit fundamental domain by the authors \cite{FL1}, \cite{FL2} for the Picard modular group $PU(2,1;{\bf Z}[i])$. This method   
uses Siegel sets on the unbounded hyperquadric model of the complex hyperbolic space. 
The spectral analysis and the scattering theory of the corresponding automorphic Laplace-Beltrami operator are  developed in \cite{FL3}.  Some of the results discussed here were announced at the conference "Geometric Analysis of PDE and Several Complex Variables" in Serra Negra, Brazil, in 2003. After completing this paper the authors received the interesting preprint of Falbel-Parker \cite{FP2}. In \cite{FP2}, a fundamental domain is  constructed for the Eisenstein-Picard modular group using very  different methods.  

The main results of the paper are the following theorems:
\begin{thm}
A fundamental domain for the Picard modular group is 
\begin{eqnarray*}
{\cal F} & \equiv & \{z\in{\bf C}^{2};\ \ 0\leq\Re e z_{1},\ 0\leq\Im m z_{1},\ \Re e z_{1}+\Im m z_{1}\leq 1,\\
\ & \ & |\Re e z_{2}|\leq\frac{1}{2},\ 
 Q_{1}(z)\equiv|z_{2}|^{2}\geq 1,\\ 
\ & \ &  Q_{3+r}(z)\equiv|r+i-(1+i)z_{1}+z_{2}|^{2}\geq 1,\ r=-1,0,1, \\
\ & \ &  Q_{5+\frac{1+r}{2}}(z)\equiv|r+i-2iz_{1}+2z_{2}|^{2}\geq 2,\ r=-1,1, \\
\ & \ &  Q_{7+\frac{1+r}{2}}(z)\equiv|r+i-2z_{1}+2z_{2}|^{2}\geq 2,\ r=-1,1 \  \}
\end{eqnarray*}
\label{maint05}
\end{thm}

Theorem \ref{maint1} contains more precise description of the structure of the fundamental domain ${\cal F}$.  

\begin{thm}

There are eight explicitly known holomorphic automorphims $G_{1}=J, G_{2},\dots, G_{8}$ in the Picard modular group such that the set 
\begin{eqnarray*}
{\cal F} & \equiv & \{z\in{\bf C}^{2};\ \ 0\leq\Re e z_{1},\ 0\leq\Im m z_{1},\ \Re e z_{1}+\Im m z_{1}\leq 1,\\
\ & \ & |\Re e z_{2}|\leq\frac{1}{2},\ 
 |z_{2}|^{2}\geq 1,\\ 
\ & \ &  |\det G'_{j}(z)|^{2}\leq 1,\ j=2,\dots, 8 \}
\end{eqnarray*}
is a fundamental domain of the Picard modular group acting on the complex hyperbolic space ${\bf CH}^{2}$. All eight transformations are needed. The transformations $G_{1},\dots,G_{8}$ can be described as follows: 

There are four transformations with dilation parameter $1$:
\begin{equation}
G_{1}(z_{1},z_{2})\equiv J(z_{1},z_{2})=(\frac{iz_{1}}{z_{2}}, -\frac{1}{z_{2}}),
\end{equation}
\begin{equation}
G_{r+3}=J\circ P_{r+3}=J\circ N_{(1+i,r+i)}\circ M_{-1},\ \ \ r=-1,0,1. 
\label{G234}
\end{equation}

There are four transformations with dilation parameter $\sqrt{2}$:

\begin{eqnarray}
G_{5+\frac{1+r}{2}} & \equiv & N_{(1,\frac{r+i}{2})}\circ J\circ P_{5+\frac{1+r}{2}} \nonumber \\
 & = & 
N_{(1,\frac{r+i}{2})}\circ J\circ N_{(-1+ri, r+i)}\circ A_{\sqrt{2}}\circ 
M_{\frac{-1+ri}{\sqrt{2}}},\ \ \ r=-1,1, 
\label{G56}
\end{eqnarray}
and 
\begin{eqnarray}
G_{7+\frac{1+r}{2}} & \equiv & N_{(i,\frac{r+i}{2})}\circ J\circ P_{7+\frac{1+r}{2}} \nonumber \\
 & = & 
N_{(i,\frac{r+i}{2})}\circ J\circ N_{(-r-i, r+i)}\circ A_{\sqrt{2}}\circ 
M_{\frac{-1+ri}{\sqrt{2}}},\ \ \ r=-1,1.
\label{G78}
\end{eqnarray}
\label{maint1}
\end{thm}
The precise definition of the holomorphic automorphisms $P$, $J$, $N$, $A$, and $M$ is described in Section 2. We mention that the functions $Q_{j}$ in Theorem \ref{maint05} are the Jacobian determinant of the transformations $G_{j}$ to some power, see also (\ref{ejacobi}) and (\ref{ejacobi2}) for more details.

Three important geometric properties of the fundamental domain ${\cal F}$ are  stated in the following theorem. These properties play a crucial role in our approach to the Laplace-Beltrami operator associated to the Picard modular group. The spectral and scattering analysis of the corresponding automorphic Laplace-Beltrami operator is developed in \cite{FL3}.

\begin{thm} (i) The fundamental domain ${\cal F}$ is invariant under the involutive transformation $S(z_{1},z_{2})=(i\bar{z_{1}}, -\bar{z_{2}})$. 

(ii) The two dimensional edge of the fundamental domain ${\cal F}$ at $z_{1}=0$ is identical to the standard fundamental domain for the modular group. More precisely,
$${\cal F}\cap \{z_{1}=0\}=\{(0,z_{2})\in{\bf C}^{2};\ |\Re ez_{2}|\leq \frac{1}{2},\ 
\Im m z_{2}>0,\  |z_{2}|\geq 1\}.$$ 

(iii) The fundamental domain $\cal F$ has a product structure near infinity, that is, for large enough $a>0$ 
\begin{eqnarray}
{\cal F}\cap\{z\in{\bf C}^{2};\ \Im m z_{2}\geq a\} & = & \{z\in{\bf C}^{2};\  \ 0\leq\Re e z_{1},\ 0\leq\Im m z_{1},
\nonumber\\
 \ & \ &  \Re e z_{1}+\Im m z_{1}\leq 1,\ |\Re e z_{2}|\leq\frac{1}{2},\ \Im m z_{2}\geq a\}.
\label{fa}
\end{eqnarray}

\label{maint3}
\end{thm}

\begin{rem}
The existence of a symmetry like $S$ in Theorem \ref{maint3} for a general discrete automorphism group has important consequences for the spectrum of the automorphic Laplace-Beltrami operator. In particular, it implies the existence of infinitely many embedded eigenvalues in the continuous spectrum. 
These results will appear elsewhere. 

\label{r2}
\end{rem}

\begin{rem}
We mention that the following identities hold for the transformations $G_{j}$, $j=1,\dots,8$: 
\begin{eqnarray}
SG_{1}S & = & G_{1},\nonumber \\ 
SG_{r+3}S & = & G_{-r+3},\ \ r=-1,0,1,\nonumber \\ SG_{5+\frac{1+r}{2}}S & = & G_{7+\frac{1-r}{2}},\ \ r=-1,1.
\label{eidS}
\end{eqnarray} 
\label{r25}
\end{rem}

\section{Preliminaries}

In this section we review the necessary background material on complex hyperbolic space and Picard modular group. A more extensive treatment of these topics may be found in \cite {E}, \cite{G}, \cite{C} and \cite{T}.
 
 The Hermitian symmetric space 
$\hbox{SU}(n,1)/\hbox{S}(\hbox{U}(1)\times \hbox{U}(n))$ is called as the complex hyperbolic space, ${\bf CH}^{n}$. 
A standard model of the complex hyperbolic space is the complex unit ball $B^{n}=\{z\in {\bf C}^{n};\ |z|<1\}$ with the Bergman metric $g=\sum_{j,k=1}^{n}g_{j,k}(z)dz_{j}\otimes d\bar{z}_{k}$, where $g_{j,k}=
\hbox{const}\cdot\partial_{j}\overline{\partial}_{k}\log(1-|z|^{2})$. 
We will use also the unbounded hyperquadric model of the complex hyperbolic space, that is $D^{n}=\{z\in {\bf C}^{n};\ \Im m z_{n}>\frac{1}{2}\sum_{j=1}^{n-1}|z_{j}|^{2}\}$.
The biholomorphic function 
$$z\mapsto(\frac{\sqrt{2}z_{1}}{iz_{n}-1},\dots,\frac{\sqrt{2}z_{n-1}}{iz_{n}-1},\frac{iz_{n}+1}{iz_{n}-1})$$ 
maps the hyperquadric $D^{n}$  onto the unit ball $B^{n}$. 

The holomorphic automorphism group of ${\bf CH}^{n}$, $\hbox{\bf Aut}({\bf CH}^{n})$, consists of 
rational functions $g=(g_{1},\dots,g_{n}):D^{n}\mapsto D^{n}$,
$$g_{j}(z)=\frac{a_{j+1,1}+\sum_{k=2}^{n+1}a_{j+1,k}z_{k-1}} 
{a_{1,1}+\sum_{k=2}^{n+1}a_{1,k}z_{k-1}},\ \ j=1,\dots,n.$$
The automorphisms act linearly in homogeneous coordinates $(\zeta_{0},\dots,\zeta_{n})$, $z_{j}=\frac{\zeta_{j}}{\zeta_{0}}$, $j=1,\dots,n$. The corresponding matrix $A=[a_{jk}]_{j,k=1}^{n+1}$ satisfies the condition 
\begin{equation}
A^{*}CA=C,
\label{e25}
\end{equation}

where 
\begin{equation*}
C\equiv\left(\begin{array}{ccc}
0 & 0& i \\
0 & I_{n-1} & 0 \\
-i & 0 & 0
\end{array}\right)
\end{equation*}
and $I_{n-1}$ is the $(n-1)\times(n-1)$ identity matrix. The determinant of the matrix $A$ is normalized to be equal to $1$. The matrix $C$ is the matrix of the quadratic form of the defining function of $D^{n}$ written in homogeneous coordinates. 

We now describe three important classes of holomorphic automorphisms stabilizing $\infty$. 
\begin{ex} Heisenberg translations. Let $a\in\partial D^{2}$. Then the Heisenberg translation $N_{a}\in
\hbox{\bf Aut}({\bf CH}^{2})$ is defined as 
$N_{a}(z_{1},z_{2})=(z_{1}+a_{1},z_{2}+a_{2}+iz_{1}\bar{a}_{1})$. If we write $a=(a_{1},a_{2})=(\gamma, r+\frac{i}{2}|\gamma|^{2})$ with $\gamma\in{\bf C}$, $r\in{\bf R}$ then the Heisenberg translation is given by 
$N_{\gamma,r+\frac{i}{2}|\gamma|^{2}}(z_{1},z_{2})=
(z_{1}+\gamma,z_{2}+r+\frac{i}{2}|\gamma|^{2}+iz_{1}\bar{\gamma})$. The corresponding matrix representation is
\begin{equation*}
N\equiv\left(\begin{array}{ccc}
1 & 0 & 0 \\
a_{1} & 1 & 0 \\
a_{2} & i\bar{a}_{1} & 1
\end{array}\right)=
\left(\begin{array}{ccc}
1 & 0 & 0 \\
\gamma & 1 & 0 \\
r+\frac{i}{2}|\gamma|^{2} & i\bar{\gamma} & 1
\end{array}\right)
\end{equation*}
 and condition (\ref{e25}) is satisfied. Let ${\cal N}\subset \hbox{SL}(3,{\bf C})$ be the set of all Heisenberg translations. 
\label{ex5}
\end{ex}

\begin{ex} Dilations. Let $\delta>0$. Then the dilation $A_{\delta}(z)=(\delta z_{1},\delta^{2}z_{2})$ is a holomorphic automorphism of $D^{2}$. Its matrix representation is 
\begin{equation*}
A\equiv\left(\begin{array}{ccc}
\frac{1}{\delta} & 0 & 0 \\
0 & 1 & 0 \\
0 & 0 & \delta
\end{array}\right).
\end{equation*}
It is easy to verify that $A$ satisfies the condition (\ref{e25}). 
We denote the set of all dilations by ${\cal A}\subset \hbox{SL}(3,{\bf C})$. 
\label{ex6}
\end{ex}

\begin{ex} Rotations. Let $\varphi\in {\bf R}$. The rotation in the first variable by $e^{i\varphi}$, $(z_{1},z_{2})\mapsto (e^{i\varphi}z_{1},z_{2})$ is a holomorphic automorphism of $D^{2}$. There are three matrices  
\begin{equation*}
M\equiv\left(\begin{array}{ccc}
\beta & 0 & 0 \\
0 & \beta^{-2} & 0 \\
0 & 0 & \beta
\end{array}\right),
\end{equation*}
 $\beta=e^{-i\varphi/3+2\pi ik/3}$, $k=0,1,2$ corresponding to the same rotation. 
It is easy to see that $M$ satisfies the condition (\ref{e25}). The set of rotations is denoted by ${\cal M}=\{M\in \hbox{SL }(3,{\bf C});\ \beta\in{\bf C},\ |\beta|=1\}$.
\label{ex7}
\end{ex}

The next transformation is an involutive automorphism. 

\begin{ex} The involution $J(z_{1},z_{2})=(\frac{iz_{1}}{z_{2}},-\frac{1}{z_{2}})$ is a holomorphic automorphism of $D^{2}$ mapping $\infty$ into $(0,0)$. A matrix representation of $J$ is 
\begin{equation*}
J\equiv\left(\begin{array}{ccc}
0 & 0 & i \\
0 & -1 & 0 \\
-i & 0 & 0
\end{array}\right).
\end{equation*}
 Notice that $J^{2}=I$, and that every point of form $(z,i)$ is a fixed point of $J$. 
\label{ex10}
\end{ex}

Let $z$ be a boundary point of $D^{2}$, i.e. $z\in\partial D^{2}\cup\{\infty\}$ The isotropy subgroup (stabilizer subgroup) $\Gamma_{z}$ of $z$ contains all the holomorphic automorphisms that leave $z$ fixed, that is $\Gamma_{z}\equiv\{g\in\hbox{\bf Aut}({\bf CH}^{2});\ g(z)=z\}$. The isotropy subgroup of $\infty$ consists of lower triangular matrices, that is 
$$\Gamma_{\infty}=\{P\in\mbox{SL}(3,{\bf C}); P^{*}CP=C, p_{12}=p_{13}=p_{23}=0\}.$$
The stabilizer group of $\infty$, ${\cal P}\equiv\Gamma_{\infty}$ can be decomposed as the product 
$${\cal P}={\cal N}{\cal A}{\cal M}.$$
This decomposition is called Langlands decomposition. An element of the stability group $P\in{\cal P}$ can be decomposed as 
\begin{eqnarray*}
P=\left(\begin{array}{ccc}
p_{11} & 0 & 0 \\
p_{21} & p_{22} & 0 \\
p_{31} & p_{32} & p_{33}
\end{array}\right)=
NAM=
\left(\begin{array}{ccc}
\frac{\beta}{\delta} & 0 & 0 \\
 \frac{\beta\gamma}{\delta}& \beta^{-2} & 0 \\
\frac{\beta}{\delta}(r+\frac{i}{2}|\gamma|^{2}) & i\bar{\gamma}\beta^{-2} & \beta\delta
\end{array}\right).
\end{eqnarray*}

The Jacobi determinant of the automorphism 
$$G(z)=\left(\frac{g_{21}+g_{22}z_{1}+g_{23}z_{2}}{g_{11}+g_{12}z_{1}+
g_{13}z_{2}},
\frac{g_{31}+g_{32}z_{1}+g_{33}z_{2}}{g_{11}+g_{12}z_{1}+
g_{13}z_{2}}\right)\in\hbox{\bf Aut}({\bf CH}^{2})$$
is given by the formula 
\begin{equation}
\det G'(z)=(g_{11}+g_{12}z_{1}+g_{13}z_{2})^{-{3}}.
\label{ejacobi}
\end{equation}
We will use the notation $Q(G)$ for the quadratic form 
\begin{equation}
Q(G)=|g_{11}+g_{12}z_{1}+g_{13}z_{2}|^{2}.
\label{ejacobi2}
\end{equation}

An analogue of the modular group $\hbox{SL}(2, {\bf Z})$ in higher dimensional complex hyperbolic space is the Picard modular group. 
\begin{defn} 
The Picard modular group is a holomorphic automorphism group of 
${\bf CH}^{n}$ defined as
\begin{equation}
\Gamma=\{A\in \hbox{SL }(n+1, {\bf C});\ \ \ A^{*}CA=C, \ a_{jk}\in{\bf Z}[i]\}.
\label{epic}
\end{equation}
\label{pic}
\end{defn}
We recall that the condition $a_{jk}\in{\bf Z}[i]$ means that the entries of the matrix $A$ are Gaussian integers, that is $\Re e a_{jk},\ \Im m a_{jk}\in{\bf Z}$. 

The Picard modular group, $\Gamma$, is a discontinuous subgroup of 
$\hbox{\bf Aut}({\bf CH}^{n})$.
\begin{rem} We note that the Heisenberg translation $N_{\gamma,r+\frac{i}{2}|\gamma|^{2}}$  
is in the Picard modular group, $\Gamma$, if and only if $r\in {\bf Z}$, and $\gamma$ is a Gaussian integer with the property that 
\begin{equation}
\frac{|\gamma|^{2}}{2}\in{\bf Z}.
\label{e1/2}
\end{equation} 
\label{r3}
\end{rem}
We close this section by recalling the definition of fundamental domain \cite{T}. Note that the fundamental domain is a closed set in this definition. 
\begin{defn}  
A set ${\cal F}\subset D^{n}$ is a fundamental domain for the automorphism group $\Gamma\subset$\hbox{\bf Aut}${\bf (CH}^{n})$ if 

(i) for all $z\in D^{n}$ there is a $G \in\Gamma$
 such that $G(z)\in{\cal F}$;

(ii) whenever $z\in{\cal F}$ and $G(z)\in{\cal F}$ for some $G\in\Gamma\setminus\{I\}$ then $z$ lies on the boundary of ${\cal F}$.  

\label{fd}
\end{defn}

We introduce the definition of {\it Siegel set} in a metric space $(X,d)$, \cite{S}. 
Let $\Gamma$ be a discrete subgroup of the isometry group of $X$. 
\begin{defn}
A closed subset ${\cal S}\subset X$ is a Siegel set for $\Gamma$ if it satisfies the following two properties:

(i) $\Gamma\cdot {\cal S}=X$, i.e. for all $x\in X$ there is a $G\in\Gamma$ such that $G(x)\in {\cal S}$;

(ii) the set $\{G\in\Gamma;\ {\cal S}\cap G({\cal S})\not=\emptyset\}$ is finite.

\label{Sie1}
\end{defn}

Let $L\geq 0$ and let $\Delta$ be the closed triangle in ${\bf C}$ with vertices $0$, $1$ and $i$. 
\begin{defn}
We define  the closed subset $S(L)$ as  
\begin{equation}
S(L)=\{z\in{\bf C}^{2};\  z_{1}\in\Delta,\ |\Re e  z_{2}|\leq\frac{1}{2},\ \Im m z_{2}-\frac{1}{2}|z_{1}|^{2}\geq L
\}.
\label{eSie2}
\end{equation}
\label{Sie2}
\end{defn}

We proved in \cite{FL1}, \cite{FL2} that if $0<L<\frac{\sqrt{3}-1}{2}$ then the set $S(L)$ is a Siegel set for the Picard modular group $\Gamma\subset$ \hbox{\bf Aut}${\bf (CH}^{2})$.

\subsection{Outline of the proof}. 

The basic idea of the proof can be described easily. In \cite{FL1} we constructed a semi-explicit 
fundamental domain for the Picard modular group $\Gamma$. There are holomorphic automorphisms $G_{1}\equiv J$, $G_{2}$, \dots, $G_{N}$ in the Picard modular group such that the set 
\begin{eqnarray}
{\cal F} & \equiv & \{z\in{\bf C}^{2};\ \ 0\leq\Re e z_{1},\ 0\leq\Im m z_{1},\ \Re e z_{1}+\Im m z_{1}\leq 1,\nonumber\\
\ & \ & |\Re e z_{2}|\leq\frac{1}{2},\ 
 |z_{2}|^{2}\geq 1,\  |\det G'_{j}(z)|^{2}\leq 1,\ j=2,\dots, N \}
\label{sef}
\end{eqnarray}
is a fundamental domain for the Picard modular group acting on the complex hyperbolic space ${\bf CH}^{2}$. The transformations $G_{j}$ $j=2,\dots, N$ satisfy the Siegel property, but are not known explicitly. 
  
Let ${\cal F}_{1}\equiv S(L)\cap\{z\in{\bf C}^{2}; |z_{2}|\geq 1\}=S(L)\cap\{z\in{\bf C}^{2}; |\det G_{1}'(z)|^{2}\leq 1\}$. Clearly ${\cal F}\subset{\cal F}_{1}$. We will prove that if $G$ is one of the transformations $G_{j}$, $j=2,\dots, N$ in the description of ${\cal F}$ in (\ref{sef}) then either 

(1) $|\det G'(z)|\leq 1$ for all $z\in{\cal F}_{1}$;

or

(2) there is a transformation $G_{j}$, $j=2,\dots,8$ appearing in (\ref{G234}), (\ref{G56}) and (\ref{G78}) such that 
$|\det G'(z)|\leq |\det G_{j}'(z)|$ for all $z\in{\cal F}_{1}$. In either case, the transformation $G$ does not contribute to the fundamental domain ${\cal F}$.

\section{Proof of the Main Theorems}. 
We recall from \cite{FL1}, \cite{FL2} that any $G=[g_{jk}]\in\Gamma\setminus\Gamma_{\infty}$ can be decomposed as $G=NJP$. The matrix $P=[p_{jk}]\in\Gamma_{\infty}$ is lower triangular and parameterized by $\delta>0$, $\beta\in{\bf C}$ with $|\beta|=1$, $r\in{\bf R}$, $\gamma\in{\bf C}$ as 
\begin{equation}
p_{31}=\frac{\beta}{\delta}(r+\frac{i}{2}|\gamma|^{2}),\ \  p_{32}=i\bar{\gamma}\beta^{-2},\ \ p_{33}=\beta\delta.
\label{e101}
\end{equation}
The transformations $N$ and $P$ in the decomposition of 
$G=[g_{jk}]\in\Gamma\setminus\Gamma_{\infty}$ are not necessarily in the Picard modular group $\Gamma$, the entries of $N$, $P$ are not necessarily Gaussian integers. However 
\begin{equation}
g_{1j}=ip_{3j},\ j=1,2,3.
\label{e1015}
\end{equation}

\begin{lem}
Furthermore, the parameters $\delta$, $r$, $\gamma$ in the decomposition of $G=[g_{jk}]\in\Gamma\setminus\Gamma_{\infty}$ satisfy the following conditions:
\begin{equation}
\delta>0,\ \ 1\leq\delta^{2}\in{\bf Z},\ \ r\in{\bf Z},\ \ |\gamma|^{2}/2\in{\bf Z},\ \ \delta^{2}\Re e\gamma\in{\bf Z},\ \ \delta^{2}\Im m\gamma\in{\bf Z}.
\label{e102}
\end{equation}

\label{l101}
\end{lem}

\noindent{\bf Proof of Lemma \ref{l101}}. Since the entries $g_{jk}$ are Gaussian integers the statement follows from the identities $g_{11}\bar{g}_{13}=r+i\frac{|\gamma|^{2}}{2}$, $1\leq |g_{13}|^{2}=\delta^{2}$, $g_{13}^{2}g_{12}=\delta^{2}\bar{\gamma}$.

We proceed to reduce the number of the transformations used to describe the fundamental domain to a handful explicitly given automorphisms. 
\subsection{Crude estimates for the automorphisms $G_{2},\dots, G_{N}$}
\begin{lem}
Let $G=[g_{jk}]\in\Gamma\setminus\Gamma_{\infty}$ be a holomorphic automorphism of ${\bf CH}^{2}$. If $|g_{13}|\geq \frac{2}{\sqrt{3}-1}$ then $|g_{11}+g_{12}z_{1}+g_{13}z_{2}|^{2}\geq 1$ for all $z\in {\cal F}_{1}$.
\label{l103}
\end{lem}

\noindent{\bf Proof of Lemma \ref{l103}}. Since $G\not\in\Gamma_{\infty}$ we have $g_{13}\not=0$ and 
\begin{eqnarray}
|g_{11}+g_{12}z_{1}+g_{13}z_{2}|^{2} & = & |g_{13}|^{2}|\Re e\frac{g_{11}}{g_{13}}+\Re e\frac{g_{12}}{g_{13}}z_{1}+\Re e z_{2}|^{2}+\nonumber\\
 & & |g_{13}|^{2}|\Im m\frac{g_{11}}{g_{13}}+\Im m
\frac{g_{12}}{g_{13}}z_{1}+\Im m z_{2}|^{2}\nonumber\\
 & \geq & |g_{13}|^{2}(\Im m\frac{g_{11}}{g_{13}}+\Im m
\frac{g_{12}}{g_{13}}z_{1}+\Im m z_{2})^{2}.
\label{e103}
\end{eqnarray}
Using that 
$$|\Im m\frac{g_{12}}{g_{13}}z_{1}|\leq \frac{1}{2}(|\frac{g_{12}}{g_{13}}|^{2}+|z_{1}|^{2})$$ 
and (\ref{e101}), (\ref{e1015}) we obtain  
\begin{eqnarray}
\Im m\frac{g_{11}}{g_{13}}+\Im m
\frac{g_{12}}{g_{13}}z_{1}+\Im m z_{2} & \geq  & \Im m\frac{g_{11}}{g_{13}}-\frac{1}{2}|\frac{g_{12}}{g_{13}}|^{2}+\Im m z_{2}-\frac{1}{2}|z_{1}|^{2}\\
\ & = & \frac{1}{|g_{13}|^{2}}(\Im m p_{31}\bar{p_{33}}-\frac{1}{2}|p_{32}|^{2})+\Im m z_{2}-\frac{1}{2}|z_{1}|^{2}\\
\ & = & \Im m z_{2}-\frac{1}{2}|z_{1}|^{2}
\label{e104}
\end{eqnarray}
Since $z\in{\cal F}_{1}$ we have $|\Re e z_{2}|\leq 1/2$ and $|z_{1}|\leq 1$ so $|\Im m z_{2}|^{2}\geq 1-|\Re e z_{2}|^{2}\geq \frac{3}{4}$. This implies that 
\begin{eqnarray*}
\Im m z_{2}-\frac{1}{2}|z_{1}|^{2}\geq\frac{\sqrt{3}-1}{2}\ \ \forall z\in{\cal F}_{1}
\end{eqnarray*}
 and proves Lemma \ref{l103}.

\begin{cor} Let $G=[g_{jk}]\in\Gamma\setminus\Gamma_{\infty}$ be a holomorphic automorphism of ${\bf CH}^{2}$. If there is a point $z^{0}\in{\cal F}_{1}$ such that $|\det G'(z^{0})|^{2}>1$ then the possible values of $|g_{13}|$ are $1$, $\sqrt{2}$, $2$ and $\sqrt{5}$. 
\label{c103}
\end{cor}

\noindent{\bf Proof of Corollary \ref{c103}}. Since 
$$|\det G'(z)|^{2}=|g_{11}+g_{12}z_{1}+g_{13}z_{2}|^{-6}$$
Lemma \ref{l103} implies that $0<|g_{13}|^{2}<\frac{4}{4-2\sqrt{3}}\approx 7.4641$. The possible values of $|g_{13}|^{2}$ for the Gaussian integer $g_{13}$ are the positive integers less than $8$ that are sum of squares of integers. These are $1$, $2$, $4$, $5$. 

\begin{lem}
Let $G=[g_{jk}]\in\Gamma\setminus\Gamma_{\infty}$ be a holomorphic automorphism of ${\bf CH}^{2}$ with $\delta\equiv|g_{13}|=1,\sqrt{2},2,\sqrt{5}$. There is a number $\gamma_{0}(\delta)=\delta+\sqrt{(1-\sqrt{3})\delta^{2}+2\delta}$ such that if $|g_{12}|\geq \gamma_{0}(\delta)$ then  $|g_{11}+g_{12}z_{1}+g_{13}z_{2}|^{2}\geq 1$ for all $z\in {\cal F}_{1}$.
\label{l106}
\end{lem}

\noindent{\bf Proof of Lemma \ref{l106}}. We estimate the right hand side of (\ref{e103}) using that for all $z=(z_{1},z_{2})\in{\cal F}_{1}$ we have $|z_{1}|\geq 1$ and $\Im m z_{2}\geq\frac{\sqrt{3}}{2}$. So 
\begin{eqnarray}
|g_{13}|(\Im m\frac{g_{11}}{g_{13}}+\Im m
\frac{g_{12}}{g_{13}}z_{1}+\Im m z_{2}) & \geq & 
|g_{13}|(\Im m\frac{g_{11}}{g_{13}}-|
\frac{g_{12}}{g_{13}}||z_{1}|+\Im m z_{2})  \nonumber\\ 
\ & \geq & \delta(\frac{|\gamma|^{2}}{2\delta^{2}}-\frac{|\gamma|}{\delta}+\frac{\sqrt{3}}{2}) \nonumber\\
\ & = & \frac{|\gamma|^{2}}{2\delta}-|\gamma|+\frac{\sqrt{3}}{2}\delta.
\label{e106}
\end{eqnarray}
The quadratic expression $\frac{|\gamma|^{2}}{2\delta}-|\gamma|+\frac{\sqrt{3}}{2}\delta$ is at least one if 
$$|\gamma|\geq \gamma_{0}(\delta)\equiv\delta+\sqrt{(1-\sqrt{3})\delta^{2}+2\delta}.$$
Therefore $|g_{11}+g_{12}z_{1}+g_{13}z_{2}|^{2}\geq 1$ for all $z\in{\cal F}_{1}$ if $|g_{12}|\geq \gamma_{0}(\delta)$.
  
\begin{rem}
The approximate values of $\gamma_{0}(\delta)$ for $\delta=1,\sqrt{2},2,\sqrt{5}$ are as follows: $\gamma_{0}(1)=2.12603$, $\gamma_{0}(\sqrt{2})=2.58226$, $\gamma_{0}(2)=3.03528$, $\gamma_{0}(\sqrt{5})=3.13711$.
\label{r103}
\end{rem}

\begin{lem}
Let $G=[g_{jk}]\in\Gamma\setminus\Gamma_{\infty}$ be a holomorphic automorphism of ${\bf CH}^{2}$ with $|g_{13}|\equiv\delta\in\{1,\sqrt{2},2,\sqrt{5}\}$, $|g_{12}|\leq \gamma_{0}(\delta)$. There is a number $r_{0}(\delta)=\delta\gamma_{0}(\delta)+\frac{\delta^{2}}{2}+\delta\sqrt{1-\delta^{2}(1-\frac{\sqrt{3}}{2})}$ such that if $|r|\equiv|\Re e g_{11}\bar{g_{13}}|\geq r_{0}(\delta)$ then  $|g_{11}+g_{12}z_{1}+g_{13}z_{2}|^{2}\geq 1$ for all $z\in {\cal F}_{1}$.
\label{l109}
\end{lem}

\noindent{\bf Proof of Lemma \ref{l109}}. Let $z\in{\cal F}_{1}$ and $|g_{12}|\leq \gamma_{0}(\delta)$. Then we have the estimate 
\begin{eqnarray*}
\frac{|r|}{\delta} & \leq & \left|\frac{r}{\delta}-\Im m \frac{\bar{\gamma}}{\beta^{3}}z_{1}+\delta\Re e z_{2}\right|+|\gamma||z_{1}|+\delta|\Re e z_{2}| \\
\ & \leq & \left|\frac{r}{\delta}-\Im m \frac{\bar{\gamma}}{\beta^{3}}z_{1}+\delta\Re e z_{2}\right|+\gamma_{0}(\delta)+\frac{\delta}{2}.
\end{eqnarray*}

So 
$$\left|\frac{r}{\delta}-\Im m \frac{\bar{\gamma}}{\beta^{3}}z_{1}+\delta\Re e z_{2}\right|^{2}\geq (\frac{|r|}{\delta}-\gamma_{0}(\delta)-\frac{\delta}{2})^{2}$$
and therefore 
\begin{eqnarray*}
|g_{11}+g_{12}z_{1}+g_{13}z_{2}|^{2} & = & |\frac{1}{\delta}(r+i\frac{|\gamma|^{2}}{2})+
i\frac{\bar{\gamma}}{\beta^{3}}z_{1}+\delta z_{2}|^{2}\\
\ & = & (\frac{r}{\delta}-\Im m \frac{\bar{\gamma}}{\beta^{3}}z_{1}+\delta\Re e  z_{2})^{2}+\\
\ & & (\frac{|\gamma|^{2}}{2\delta}+\Re e  \frac{\bar{\gamma}}{\beta^{3}}z_{1}+\delta\Im m  z_{2})^{2}\\
\ & \geq & (\frac{|r|}{\delta}-\gamma_{0}(\delta)-\frac{\delta}{2})^{2}+(\delta\frac{\sqrt{3}-1}{2})^{2}
\end{eqnarray*}
according to the proof of Lemma \ref{l103}. The right hand side is at least $1$ for all $z\in{\cal F}$ provided $|r|\geq r_{0}(\delta)$.

\begin{rem}
The approximate values of $r_{0}(\delta)$ for $\delta=1,\sqrt{2},2,\sqrt{5}$ are as follows: $r_{0}(1)=3.55664$, $r_{0}(\sqrt{2})=5.86186$, $r_{0}(2)=9.43305$, $r_{0}(\sqrt{5})=10.7996$.
\label{r106}
\end{rem}

\begin{lem}
Let $G=[g_{jk}]\in\Gamma\setminus\Gamma_{\infty}$ be a holomorphic automorphism of ${\bf CH}^{2}$ with $g_{12}=ip_{32}=-\bar{\gamma}\frac{1}{\beta^{2}}=0$, that is $\gamma=0$. Then  $|g_{11}+g_{12}z_{1}+g_{13}z_{2}|^{2}\geq 1$ for all $z\in {\cal F}_{1}$.
\label{l112}
\end{lem}

\noindent{\bf Proof of Lemma \ref{l112}}. If $\gamma=0$ then 
\begin{eqnarray*}
|g_{11}+g_{12}z_{1}+g_{13}z_{2}|^{2} & = & |\frac{\beta}{\delta}r+\beta\delta z_{2}|^{2}\\
 & \geq & \frac{r^{2}}{\delta^{2}}-2|r||\Re e z_{2}|+\delta^{2}|z_{2}|^{2}\\
 & \geq & \frac{r^{2}}{\delta^{2}}-|r|+\delta^{2}. 
\end{eqnarray*}

When $r=0$ then the statement follows from $\delta\geq 1$. If $r\in{\bf Z}\setminus\{0\}$ then $|r|\geq 1$ so $\frac{r^{2}}{\delta^{2}}-|r|+\delta^{2}=(\frac{|r|}{\delta}-\delta)^{2}+|r|\geq 1$. 

The next lemma rules out the case $\delta=\sqrt{5}$.
\begin{lem}
Let $G=[g_{jk}]\in\Gamma\setminus\Gamma_{\infty}$ be a holomorphic automorphism of ${\bf CH}^{2}$ with $|g_{13}|=|ip_{33}|=\delta=\sqrt{5}$. Then  $|g_{11}+g_{12}z_{1}+g_{13}z_{2}|^{2}\geq 1$ for all $z\in {\cal F}_{1}$.
\label{l115}
\end{lem}

\noindent{\bf Proof of Lemma \ref{l115}}. According to Remark \ref{r103} and Lemma \ref{l112} the statement of Lemma \ref{l115} is valid when $g_{12}=0$ or $|g_{12}|\geq\gamma_{0}(\sqrt{5})\approx 3.13711$. So it is enough to consider the case when $0<|g_{12}|^{2}<\gamma_{0}(\sqrt{5})^{2}<10$. Since $g_{11}\bar{g_{13}}=r+i\frac{|\gamma|^{2}}{2}$ is a Gaussian integer so $|g_{12}|^{2}=|\gamma|^{2}\in 2{\bf Z}$. Then it is easy to see that the possible values of $|g_{12}|^{2}=|\gamma|^{2}$ are $2$, $4$, $8$. Note that 
$|\bar{g}_{13}g_{12}|^{2}=\delta^{2}|\gamma|^{2}=
5|\gamma|^{2}$, so $|\Re e (\bar{g}_{13}g_{12})|^{2}+|\Im m (\bar{g}_{13}g_{12})|^{2}=10,20,40$ with $(|\Re e (\bar{g}_{13}g_{12})|, |\Im m (\bar{g}_{13}g_{12})|)\in{\bf Z}\times{\bf Z}$. One can determine the possible values of the pair of integers $(|\Re e (\bar{g}_{13}g_{12})|, |\Im m (\bar{g}_{13}g_{12})|)$ easily: $(1,3)$, $(3,1)$, $(2,4)$, $(4,2)$, $(2,6)$, and $(6,2)$. We consider these six cases separately in estimating equation (\ref{e103}). Using the notations $z_{1}=x_{1}+iy_{1}$, $z_{2}=x_{2}+iy_{2}$ we have 
\begin{eqnarray*}
\Im m\frac{g_{11}}{g_{13}}+\Im m
\frac{g_{12}}{g_{13}}z_{1}+\Im m z_{2} & = & \frac{|\gamma|^{2}}{10}+\frac{1}{5}\Im m(g_{12}\bar{g}_{13})x_{1}+\frac{1}{5}\Re e(g_{12}\bar{g}_{13})y_{1}+y_{2}\\
 & \geq & \frac{|\gamma|^{2}}{10}-\frac{1}{5}|\Im m(g_{12}\bar{g}_{13})|x_{1}-\frac{1}{5}|\Re e(g_{12}\bar{g}_{13})|y_{1}+y_{2}.
\end{eqnarray*}  

CASE 1. If $(|\Re e (\bar{g}_{13}g_{12})|, |\Im m (\bar{g}_{13}g_{12})|)=(1,3)$ then $|\gamma|^{2}=2$ so 
$$\frac{|\gamma|^{2}}{10}-\frac{1}{5}|\Im m(g_{12}\bar{g}_{13})|x_{1}-\frac{1}{5}|\Re e(g_{12}\bar{g}_{13})|y_{1}+y_{2}=\frac{1}{5}-\frac{3}{5}x_{1}-\frac{1}{5}y_{1}+y_{2}.$$
Since $z_{1}\in\Delta$, $\Im m z_{2}\geq\frac{\sqrt{3}}{2}$
we have 
$$\frac{1}{5}-\frac{3}{5}x_{1}-\frac{1}{5}y_{1}+y_{2}=\frac{1}{5}-\frac{1}{5}x_{1}-\frac{1}{5}y_{1}-\frac{2}{5}x_{1}+y_{2}\geq -\frac{2}{5}+\frac{\sqrt{3}}{2}>0.$$
Therefore
$$|g_{11}+g_{12}z_{1}+g_{13}z_{2}|^{2}\geq 5(\frac{\sqrt{3}}{2}-\frac{2}{5})^{2}\geq 1$$
 for all $z\in {\cal F}_{1}$.

CASE 2. The case $(|\Re e (\bar{g}_{13}g_{12})|, |\Im m (\bar{g}_{13}g_{12})|)=(3,1)$ is completely similar to the previous case. 

CASE 3. If $(|\Re e (\bar{g}_{13}g_{12})|, |\Im m (\bar{g}_{13}g_{12})|)=(2,4)$ then $|\gamma|^{2}=4$ so 
$$\frac{|\gamma|^{2}}{10}-\frac{1}{5}|\Im m(g_{12}\bar{g}_{13})|x_{1}-\frac{1}{5}|\Re e(g_{12}\bar{g}_{13})|y_{1}+y_{2}=\frac{2}{5}-\frac{4}{5}x_{1}-\frac{2}{5}y_{1}+y_{2}.$$
Since $z_{1}\in\Delta$, $\Im m z_{2}\geq\frac{\sqrt{3}}{2}$
we have 
$$\frac{2}{5}-\frac{4}{5}x_{1}-\frac{2}{5}y_{1}+y_{2}=\frac{2}{5}-\frac{2}{5}x_{1}-\frac{2}{5}y_{1}-\frac{2}{5}x_{1}+y_{2}\geq -\frac{2}{5}+\frac{\sqrt{3}}{2}>0.$$
Therefore
$$|g_{11}+g_{12}z_{1}+g_{13}z_{2}|^{2}\geq 5(\frac{\sqrt{3}}{2}-\frac{2}{5})^{2}\geq 1$$
 for all $z\in {\cal F}_{1}$.

CASE 4. The case $(|\Re e (\bar{g}_{13}g_{12})|, |\Im m (\bar{g}_{13}g_{12})|)=(4,2)$ is completely similar to the previous case. 

CASE 5. If $(|\Re e (\bar{g}_{13}g_{12})|, |\Im m (\bar{g}_{13}g_{12})|)=(2,6)$ then $|\gamma|^{2}=8$ so 
$$\frac{|\gamma|^{2}}{10}-\frac{1}{5}|\Im m(g_{12}\bar{g}_{13})|x_{1}-\frac{1}{5}|\Re e(g_{12}\bar{g}_{13})|y_{1}+y_{2}=\frac{4}{5}-\frac{6}{5}x_{1}-\frac{2}{5}y_{1}+y_{2}.$$
Since $z_{1}\in\Delta$, $\Im m z_{2}\geq\frac{\sqrt{3}}{2}$
we have 
$$\frac{4}{5}-\frac{6}{5}x_{1}-\frac{2}{5}y_{1}+y_{2}=\frac{4}{5}-\frac{4}{5}x_{1}-\frac{4}{5}y_{1}+\frac{2}{5}y_{1}-\frac{2}{5}x_{1}+y_{2}\geq -\frac{2}{5}+\frac{\sqrt{3}}{2}>0.$$
Therefore
$$|g_{11}+g_{12}z_{1}+g_{13}z_{2}|^{2}\geq 5(\frac{\sqrt{3}}{2}-\frac{2}{5})^{2}\geq 1$$
 for all $z\in {\cal F}_{1}$.

CASE 6. The case $(|\Re e (\bar{g}_{13}g_{12})|, |\Im m (\bar{g}_{13}g_{12})|)=(6,2)$ is completely similar to the previous case. 

This completes the proof of Lemma \ref{l115}.

\subsection{Fine estimates for the automorphisms $G_{2},\dots, G_{N}$}

The goal of this subsection to obtain finer restrictions on Jacobian determinants of the possible automorphisms contributing to the unknown part of the fundamental domain. 
There are three different classes of possible automorphisms 
since the dilation parameter of these automorphisms can have only three values, $1$, $\sqrt{2}$, $2$. In the next three subsections we analyze these three classes.
The identity
\begin{eqnarray}
Q(G)\equiv|g_{11}+g_{12}z_{1}+g_{13}z_{2}|^{2} & = & |\frac{1}{\delta}(r+i\frac{|\gamma|^{2}}{2})+i\bar{\gamma} \beta^{-3}z_{1}+\delta z_{2}|^{2}\nonumber\\
 & = & (\frac{r}{\delta} -\Im m(\frac{\bar{\gamma}}{\beta^{3}}z_{1})+\delta 
\Re e z_{2})^{2}+\nonumber\\
 &  & (\frac{|\gamma|^{2}}{2\delta}+\Re e (\frac{\bar{\gamma}}{\beta^{3}}z_{1})+\delta\Im m z_{2})^{2}\nonumber\\
& = & \left(\frac{r}{\delta} -\Im m(\frac{\bar{\gamma}}{\beta^{3}})x_{1}-\Re e (\frac{\bar{\gamma}}{\beta^{3}})y_{1}+\delta x_{2}\right)^{2}+\nonumber\\
 &  & \left(\frac{|\gamma|^{2}}{2\delta}+\Re e (\frac{\bar{\gamma}}{\beta^{3}})x_{1}-\Im m(\frac{\bar{\gamma}}{\beta^{3}})y_{1}+\delta y_{2}\right)^{2},
\label{e109}
\end{eqnarray}
for the Jacobi determinant of a transformation $G=[g_{jk}]$, which is based on (\ref{e101}), (\ref{e1015}) will be used several times.  

\subsubsection{Class of automorphisms with dilation $\delta=1$}

In this section we determine all the transformations in the Picard modular group with dilation $\delta=1$ which contribute to the boundary of the fundamental domain. 
It follows from Remark \ref{r103}, Lemma \ref{l112} and (\ref{e102}) that $0<|\gamma|^{2}<5$ and therefore the possible values of $|\gamma|^{2}$ are $2$ and $4$. Moreover $\delta=1$ implies that $\beta=\pm 1,\pm i$.  

In the next lemma we eliminate the case $|\gamma|^{2}=4$. 

\begin{lem}
Let $G=[g_{jk}]\in\Gamma\setminus\Gamma_{\infty}$ be a holomorphic automorphism of ${\bf CH}^{2}$ with $|g_{13}|=|ip_{33}|=\delta=1$ and $|g_{12}|^{2}=|\gamma|^{2}=4$. Then  $|g_{11}+g_{12}z_{1}+g_{13}z_{2}|^{2}\geq 1$ for all $z\in {\cal F}_{1}$.
\label{l116}
\end{lem}

\noindent{\bf Proof of Lemma \ref{l116}}. Since $g_{13}=i\delta\beta$ with $|g_{13}|=\delta=1$ then $\beta=\pm 1,\pm i$. It follows from $g_{12}=-\bar{\gamma}\beta^{-2}$ that $\gamma$ is a Gaussian integer, thus $\gamma=\pm 2, \pm 2i$. So the values of $\bar{\gamma}\beta^{-3}$ are $\pm 2,\pm 2i$. If $\bar{\gamma}\beta^{-3}=2$ then 
$$Q(G)=(r-2y_{1}+x_{2})^{2}+(2+2x_{1}+y_{2})^{2}\geq (2+2x_{1}+y_{2})^{2}\geq (2+\frac{\sqrt{3}}{2})^{2}\geq 1$$ 
for all $z\in{\cal F}_{1}$. Similarly, if $\bar{\gamma}\beta^{-3}=-2i$ then 
$$Q(G)=(r+2x_{1}+x_{2})^{2}+(2+2y_{1}+y_{2})^{2}\geq (2+2y_{1}+y_{2})^{2}\geq (2+\frac{\sqrt{3}}{2})^{2}\geq 1$$ 
for all $z\in{\cal F}_{1}$. It is easy to see by interchanging $x_{1}$ and $y_{1}$, $r$ and $-r$, $x_{2}$ and $-x_{2}$ that the remaining cases $\bar{\gamma}\beta^{-3}=-2, 2i$ can be estimated similarly. We will consider the case $\bar{\gamma}\beta^{-3}=-2$ only. In this case we have to estimate 
$$Q(G)=(r+2y_{1}+x_{2})^{2}+(2-2x_{1}+y_{2})^{2}.$$
If $r\geq 1$ then $r+2y_{1}+x_{2}\geq 1+2y_{1}+x_{2}\geq 1+x_{2}\geq \frac{1}{2}$ so $(r+2y_{1}+x_{2})^{2}\geq \frac{1}{4}$. Moreover $2-2x_{1}+y_{2}\geq y_{2}\geq \frac{\sqrt{3}}{2}$ so $2-2x_{1}+y_{2}\geq \frac{3}{4}$ and therefore $Q(G)\geq 1$ for all $z\in{\cal F}_{1}$. 
If $r\leq 0$ then we first estimate $Q(G)$ in the region $-\frac{1}{2}\leq x_{2}<0$. Since $2-2x_{1}\geq 0$ and $y_{2}\geq\sqrt{1-x_{2}^{2}}$ for all $z\in{\cal F}_{1}$ we 
have 
\begin{eqnarray*}
Q(G) & \geq & (r+2y_{1}+x_{2})^{2}+(2-2x_{1}+\sqrt{1-x_{2}^{2}})^{2}\\
 & = & (r+2y_{1})^{2}+2(r+2y_{1})x_{2}+x_{2}^{2}+\\
 &  & (2-2x_{1})^{2}+2(2-2x_{1})\sqrt{1-x_{2}^{2}}+1-x_{2}^{2}\\
 & \geq & 1+2(r+2y_{1})x_{2}+2(2-2x_{1})\sqrt{1-x_{2}^{2}}\\
 & \geq & 1+4y_{1}x_{2}+4(1-x_{1})\sqrt{1-x_{2}^{2}}\\
 & = & 1+4y_{1}(x_{2}+\sqrt{1-x_{2}^{2}})+4(1-x_{1}-y_{1})\sqrt{1-x_{2}^{2}}
\end{eqnarray*}
because $rx_{2}\geq 0$ in ${\cal F}_{1}$. The last expression is nonnegative since $1-x_{1}-y_{1}\geq 0$, $y_{1}\geq 0$ and $x_{2}+\sqrt{1-x_{2}^{2}}\geq -\frac{1}{2}+\sqrt{1-(1/2)^{2}}=\frac{\sqrt{3}-1}{2}>0$ in  ${\cal F}_{1}$. The same argument can be applied for $r=0$, $|x_{2}|\leq 1/2$. Now we consider $Q(G)$ in the region $0\leq x_{2}\leq \frac{1}{2}$, $r\leq -1$. It is easy to estimate $Q(G)$ if $0\leq x_{1}\leq \frac{2+\sqrt{3}}{4}$. Indeed, then 
$2-2x_{1}+y_{2}\geq 2-\frac{2+\sqrt{3}}{2}+\frac{\sqrt{3}}{2}=1$ so $Q(G)\geq (2-2x_{1}+y_{2})^{2}\geq 1$ in ${\cal F}_{1}$. If $\frac{2+\sqrt{3}}{4}\leq x_{1}\leq 1$ then using $y_{1}\geq 1-x_{1}$ we have $r+2y_{1}+x_{2}\leq -1+2-2x_{1}+\frac{1}{2}=\frac{3}{2}-2x_{1}\leq \frac{1-\sqrt{3}}{2}<0$. So we get the lower bounds $(r+2y_{1}+x_{2})^{2}\geq 
(\frac{3}{2}-2x_{1})^{2}$ and $Q(G)\geq  
(\frac{3}{2}-2x_{1})^{2}+(2-2x_{1}+y_{2})^{2}\geq (\frac{3}{2}-2x_{1})^{2}+(2-2x_{1}+\frac{\sqrt{3}}{2})^{2}$. The function on the right-hand side $f(x_{1})$ is strictly decreasing on the interval $[\frac{2+\sqrt{3}}{4},1]$, because 
$$f'(x_{1})=4\left(4x_{1}-\frac{7+\sqrt{3}}{2}\right)<0$$
provided $x_{1}<\frac{7+\sqrt{3}}{8}=1.09151...$. Therefore $Q(G)\geq f(1)=1$ for all $z\in{\cal F}_{1}$, $r\leq -1$, $0\leq x_{2}\leq \frac{1}{2}$. 
This completes the proof of Lemma \ref{l116}.

Consider case $|\gamma|^{2}=2$.  Our next step is to improve Lemma \ref{l109}.

\begin{lem}
Let $G=[g_{jk}]\in\Gamma\setminus\Gamma_{\infty}$ be a holomorphic automorphism of ${\bf CH}^{2}$ with $|g_{13}|\equiv\delta=1$, $|g_{12}|^{2}\equiv|\gamma|^{2}=2$. If $|r|\equiv|\Re e g_{11}\bar{g_{13}}|>1$ then  $|g_{11}+g_{12}z_{1}+g_{13}z_{2}|^{2}\geq 1$ for all $z\in {\cal F}_{1}$.
\label{l117}
\end{lem}

\noindent{\bf Proof of Lemma \ref{l117}}. In this case, it follows from (\ref{e109}) that we have to estimate 
\begin{eqnarray}
Q(G) & = & \left(r -\Im m(\frac{\bar{\gamma}}{\beta^{3}})x_{1}-\Re e (\frac{\bar{\gamma}}{\beta^{3}})y_{1}+ x_{2}\right)^{2}+\nonumber\\
 &  & \left(1+\Re e (\frac{\bar{\gamma}}{\beta^{3}})x_{1}-\Im m(\frac{\bar{\gamma}}{\beta^{3}})y_{1}+ y_{2}\right)^{2}. 
\label{e110}
\end{eqnarray}
It follows from (\ref{e102}) that $\gamma$ is a Gaussian integer, so $\gamma=1+i$, $1-i$, $-1+i$, $-1-i$. Since $\beta=\pm 1, \pm i$ we have $\bar{\gamma}\beta^{-3}=1+i$, $1-i$, $-1+i$, $-1-i$. 
Therefore 
$$2\leq|r|\leq |r -\Im m(\frac{\bar{\gamma}}{\beta^{3}})x_{1}-\Re e (\frac{\bar{\gamma}}{\beta^{3}})y_{1}+ x_{2}|+|-\Im m(\frac{\bar{\gamma}}{\beta^{3}})x_{1}-\Re e (\frac{\bar{\gamma}}{\beta^{3}})y_{1}+ x_{2}|\leq $$
$$|r -\Im m(\frac{\bar{\gamma}}{\beta^{3}})x_{1}-\Re e (\frac{\bar{\gamma}}{\beta^{3}})y_{1}+ x_{2}|+x_{1}+y_{1}+|x_{2}|\leq |r -\Im m(\frac{\bar{\gamma}}{\beta^{3}})x_{1}-\Re e (\frac{\bar{\gamma}}{\beta^{3}})y_{1}+ x_{2}|+\frac{3}{2},$$
so 
$$\frac{1}{2}\leq |r -\Im m(\frac{\bar{\gamma}}{\beta^{3}})x_{1}-\Re e (\frac{\bar{\gamma}}{\beta^{3}})y_{1}+ x_{2}|.$$ 
Moreover
$$1+\Re e (\frac{\bar{\gamma}}{\beta^{3}})x_{1}-\Im m(\frac{\bar{\gamma}}{\beta^{3}})y_{1}+ y_{2}\geq 1-x_{1}-y_{1}+y_{2}\geq y_{2}\geq \frac{\sqrt{3}}{2}$$
in ${\cal F}_{1}$. Therefore 
$$Q(G)\geq \left(\frac{1}{2}\right)^{2}+\left(\frac{\sqrt{3}}{2}\right)^{2}=
1$$ 
for all $z\in{\cal F}_{1}$. This proves Lemma \ref{l117}. 

\begin{lem}
Let $G=[g_{jk}]\in\Gamma\setminus\Gamma_{\infty}$ be a holomorphic automorphism of ${\bf CH}^{2}$ with $|g_{13}|\equiv\delta=1$, $|g_{12}|^{2}\equiv|\gamma|^{2}=2$. If $\bar{\gamma}\beta^{-3}=1+i, -1-i, 1-i$, and  
$r\equiv\Re e g_{11}\bar{g_{13}}=-1,0,1$ then  $|g_{11}+g_{12}z_{1}+g_{13}z_{2}|^{2}\geq 1$ for all $z\in {\cal F}_{1}$.
\label{l1175}
\end{lem}

\noindent{\bf Proof of Lemma \ref{l1175}}. If $\bar{\gamma}\beta^{-3}=1-i$ then 
\begin{eqnarray*}
Q(G) & = & (r+x_{1}-y_{1}+x_{2})^{2}+(1+x_{1}+y_{1}+y_{2})^{2}\\
 & \geq & (1+x_{1}+y_{1}+y_{2})^{2}\geq (1+y_{2})^{2}\geq \left(1+\frac{\sqrt{3}}{2}\right)^{2}\geq 1.
\end{eqnarray*}
The other two cases $\bar{\gamma}\beta^{-3}=-1-i$, and $1+i$ are similar, we will consider $\bar{\gamma}\beta^{-3}=-1-i$ only. In this case 
$$Q(G)  =  (r+x_{1}+y_{1}+x_{2})^{2}+
(1-x_{1}+y_{1}+y_{2})^{2}.$$
We may assume that $\frac{\sqrt{3}}{2}\leq x_{1}\leq 1$ because $1-x_{1}+y_{1}+y_{2}\geq 1-x_{1}+y_{2}\geq 1-x_{1}+\frac{\sqrt{3}}{2}\geq 1$ and $Q(G)\geq (1-x_{1}+y_{1}+y_{2})^{2}\geq 1$ otherwise. 
If $r=-1$ then 
\begin{eqnarray*}
Q(G) & = & (-1+x_{1}+y_{1}+x_{2})^{2}+
(1-x_{1}+y_{1}+y_{2})^{2}\\
 & = & (-1+x_{1}+y_{1})^{2}+2(-1+x_{1}+y_{1})x_{2}+x_{2}^{2}+\\
 &  & (1-x_{1}+y_{1})^{2}+2(1-x_{1}+y_{1})y_{2}+y_{2}^{2}\\
 & \geq & 2(-1+x_{1}+y_{1})x_{2}+2(1-x_{1}+y_{1})y_{2}+1\\
 & \geq & 2(1-x_{1}-y_{1})(-x_{2}+y_{2})+4y_{1}y_{2}+1\\
 & \geq & 2(1-x_{1}-y_{1})
(-\frac{1}{2}+\frac{\sqrt{3}}{2})+1\\
 & \geq & 1
\end{eqnarray*}
because $1-x_{1}-y_{1}\geq 0$, $x_{2}\geq -\frac{1}{2}$, $y_{2}\geq \frac{3}{2}$ in ${\cal F}_{1}$. 
If $r=0,1$ then $r+x_{1}+y_{1}+x_{2}\geq x_{1}+x_{2}\geq x_{1}-\frac{1}{2}\geq \frac{\sqrt{3}-1}{2}>0$ so 
$(r+x_{1}+y_{1}+x_{2})^{2}\geq (x_{1}-\frac{1}{2})^{2}$. Moreover 
$1-x_{1}+y_{1}+y_{2}\geq 1-x_{1}+y_{2}\geq 1-x_{1}+\frac{\sqrt{3}}{2}\geq \frac{\sqrt{3}}{2}$ so 
$(1-x_{1}+y_{1}+y_{2})^{2}\geq (1+\frac{\sqrt{3}}{2}-x_{1})^{2}$. Therefore 
$Q(G)\geq (x_{1}-\frac{1}{2})^{2}+
(1+\frac{\sqrt{3}}{2}-x_{1})^{2}\equiv P(x_{1})$. The minimum of the polynomial $P(x_{1})$ occurs at $x_{1}=\frac{3+\sqrt{3}}{4}> 1$. Therefore $Q(G)\geq P(x_{1})\geq P(1)=1$ for $z\in{\cal F}_{1}$. This completes the proof of Lemma \ref{l1175}.

The remaining case is $\bar{\gamma}\beta^{-3}=-1+i$. In this case $g_{11}=i\beta(r+i)$, $g_{12}=-\beta(-1+i)$, $g_{13}=i\beta$ with $r=-1,0,1$ and $\beta=\pm 1,\pm i$. So 
$Q(G)=|g_{11}+g_{12}z_{1}+g_{13}z_{2}|^{2}=|r+i-(1+i)z_{1}+z_{2}|^{2}$, and therefore 
$$|\det G'(z)|^{2}=|r+i-(1+i)z_{1}+z_{2}|^{-6},\ \ r=-1,0,1.$$

The stability component of $G$ can be determined from the Langlands decomposition:
$$P_{3+r}=\left(\begin{array}{ccc}
\beta & 0 & 0 \\
-(1+i)\beta^{2} & \beta^{-2} & 0 \\ 
\beta(r+i) &  -(1+i)\beta & \beta 
\end{array}\right)$$
with $r=-1,0,1$. Select $\beta=-1$. Consider the transformations $G_{3+r}=JP_{3+r}$. Then the Jacobian determinant of $G_{3+r}$ is 
$$|\det G_{3+r}'(z)|^{2}=|r+i-(1+i)z_{1}+z_{2}|^{-6},\ \ r=-1,0,1.$$
The transformation $G_{3+r}$ satisfies the Siegel property, namely there are points $R_{3+r}\in S(L)$ such that $G_{3+r}(R_{3+r})\in S(L)$ for $r=-1,0,1$. Indeed, it is easy to see that one can select $R_{2}=(\frac{7}{8}+i\frac{1}{16},-\frac{3}{16}+i(1+\epsilon))$,  $R_{3}=(\frac{1}{4}(1+i),i(1+\epsilon))$, $R_{4}=(\frac{1}{16}+i\frac{7}{8},\frac{3}{16}+i(1+\epsilon))$ for sufficiently small $\epsilon>0$.

\subsubsection{Class of automorphisms with dilation $\delta=\sqrt{2}$}

Let $G=[g_{jk}]\in\Gamma\setminus\Gamma_{\infty}$ be a holomorphic automorphism of ${\bf CH}^{2}$ with dilation parameter  $\delta=|g_{13}|=\sqrt{2}$. Then it follows from 
(\ref{e101}) and (\ref{e1015}) that $|p_{33}|^{2}=(\Re e p_{33})^{2}+(\Im m p_{33})^{2}=\delta^{2}=2$, so $\Re e p_{33}=\pm 1$, $\Im m p_{33}=\pm 1$, that is the values of $p_{33}$ are $p_{33}=-1-i,\ -1+i,\ 1-i,\ 1+i$. This determines the values of $\beta$:
\begin{equation}
\beta=\frac{-1-i}{\sqrt{2}},\ \frac{-1+i}{\sqrt{2}},\ \frac{1-i}{\sqrt{2}},\ \frac{1+i}{\sqrt{2}}.
\label{2beta}
\end{equation}

Note that the corresponding values of $\beta^{2}$ are $i, -i, -i, i$. Using (\ref{e101}) again we get from $\bar{\gamma}=-i\beta^{2}p_{32}$ that $\gamma$ is a Gaussian integer. Since $|\gamma|^{2}=2\Im m (p_{31}\bar{p}_{33})\in 2{\bf Z}$ and $0<|\gamma|^{2}<\gamma_{0}(\sqrt{2})^{2}<2.6^{2}<7$ according to Remark \ref{r103} and Lemma \ref{l112}, the possible values of $|\gamma|^{2}$ are $2$ and $4$. Indeed $6$ is not a possible value since it can not be written as sum of squares of two integers.   

Our next step is to eliminate the value $4$ for $|\gamma|^{2}$. 

\begin{lem}
Let $G=[g_{jk}]\in\Gamma\setminus\Gamma_{\infty}$ be a holomorphic automorphism of ${\bf CH}^{2}$ with $|g_{13}|=|ip_{33}|=\delta=\sqrt{2}$ and $|g_{12}|^{2}=|\gamma|^{2}=4$. Then  $|g_{11}+g_{12}z_{1}+g_{13}z_{2}|^{2}\geq 1$ for all $z\in {\cal F}_{1}$.
\label{l118}
\end{lem}

\noindent{\bf Proof of Lemma \ref{l118}.} We recall that 
\begin{eqnarray}
Q(G) & \geq &  
\delta^{2}\left(\frac{\frac{|\gamma|^{2}}{2}}{\delta^{2}}+
\frac{1}{\delta^{2}}\Im m(p_{32}\bar{p}_{33})x_{1}+ 
\frac{1}{\delta^{2}}\Re e(p_{32}\bar{p}_{33})y_{1}+y_{2}\right)^{2} \nonumber \\
\ & = & 2(1+\frac{1}{2}\Im m(p_{32}\bar{p}_{33})x_{1}+\frac{1}{2}\Re e(p_{32}\bar{p}_{33})y_{1}+y_{2})^{2}. 
\label{e115}
\end{eqnarray}  
Since $|\Re e(p_{32}\bar{p}_{33})|^{2}+|\Im m(p_{32}\bar{p}_{33})|^{2}=|p_{32}\bar{p}_{33}|^{2}=
\delta^{2}|\gamma|^{2}=8$ we know that $|\Re e(p_{32}\bar{p}_{33})|=2$ and $|\Im m(p_{32}\bar{p}_{33})|=2$. This gives the lower bound 
$$1-\frac{1}{2}|\Im m(p_{32}\bar{p}_{33})|x_{1}-\frac{1}{2}|\Re e(p_{32}\bar{p}_{33})|y_{1}+y_{2}=1-x_{1}-y_{1}+y_{2}\geq y_{2}\geq \frac{\sqrt{3}}{2}$$
 for the quantity in the parentheses on the right hand side of (\ref{e115}). Therefore 
$$|g_{11}+g_{12}z_{1}+g_{13}z_{2}|^{2}\geq \delta^{2}(1-x_{1}-y_{1}+y_{2})^{2}\geq 2(\frac{\sqrt{3}}{2})^{2}>1$$
 for all $z\in {\cal F}_{1}$. This completes the proof of Lemma \ref{l118}.

Now we turn our attention to the case $\delta=\sqrt{2}$, $|\gamma|^{2}=2$. In this case 
$|\Re e(p_{32}\bar{p}_{33})|^{2}+|\Im m(p_{32}\bar{p}_{33})|^{2}=|p_{32}\bar{p}_{33}|^{2}=
\delta^{2}|\gamma|^{2}=4$. Then we know that the pair of integers 
$(|\Re e(p_{32}\bar{p}_{33})|, |\Im m(p_{32}\bar{p}_{33})|)$ is equal to $(0,2)$ or $(2,0)$ so $p_{32}\bar{p}_{33}=-2i, 2i, -2, 2$. 
The next lemma further reduces the possible values of $p_{32}\bar{p}_{33}$.

\begin{lem}
Let $G=[g_{jk}]\in\Gamma\setminus\Gamma_{\infty}$ be a holomorphic automorphism of ${\bf CH}^{2}$ with $|g_{13}|=|ip_{33}|=\delta=\sqrt{2}$, $|g_{12}|^{2}=|\gamma|^{2}=2$ and $g_{12}\bar{g}_{13}=2i, 2$. Then  $|g_{11}+g_{12}z_{1}+g_{13}z_{2}|^{2}\geq 1$ for all $z\in {\cal F}_{1}$.
\label{l121}
\end{lem}

\noindent{\bf Proof of Lemma \ref{l121}.} Similarly to (\ref{e115}) we have 
$$|g_{11}+g_{12}z_{1}+g_{13}z_{2}|^{2} \geq   
2(\frac{1}{2}+\frac{1}{2}\Im m(p_{32}\bar{p}_{33})x_{1}+\frac{1}{2}\Re e(p_{32}\bar{p}_{33})y_{1}+y_{2})^{2}$$
in this case. If $g_{12}\bar{g}_{13}=p_{32}\bar{p}_{33}=2i$ then we get the lower bound 
$$2(\frac{1}{2}+x_{1}+y_{2})^{2}\geq 2(\frac{1}{2}+\frac{\sqrt{3}}{2})^{2}\geq 1.$$
The other case, $g_{12}\bar{g}_{13}=p_{32}\bar{p}_{33}=2$,  gives the lower bound 
$$2(\frac{1}{2}+y_{1}+y_{2})^{2}\geq 2(\frac{1}{2}+\frac{\sqrt{3}}{2})^{2}\geq 1.$$
These lower bounds imply Lemma \ref{l121}. 

In the next lemma we analyze the remaining two values of $g_{12}\bar{g}_{13}=p_{32}\bar{p}_{33}$, $-2$, $-2i$.

\begin{lem}
Let $G=[g_{jk}]\in\Gamma\setminus\Gamma_{\infty}$ be a holomorphic automorphism of ${\bf CH}^{2}$ with $|g_{13}|=|ip_{33}|=\delta=\sqrt{2}$, $|g_{12}|^{2}=|\gamma|^{2}=2$, $g_{12}\bar{g}_{13}=-2i, -2$.

(i) Then $r=\Re e (g_{11}\bar{g}_{13})\in{\bf Z}$ is an odd integer.

(ii) If $|r|=|\Re e(g_{11}\bar{g}_{13})|\geq 3$ then  $|g_{11}+g_{12}z_{1}+g_{13}z_{2}|^{2}\geq 1$ for all $z\in {\cal F}_{1}$.

(iii) The possible values of $r=\Re e (g_{11}\bar{g}_{13})\in{\bf Z}$ are $\pm 1$.  
\label{l124}
\end{lem}

\noindent{\bf Proof of Lemma \ref{l124}.} To prove (i) we
start with the formula $p_{31}=\frac{\beta}{\delta}(r+i\frac{|\gamma|^{2}}{2})=\frac{\beta}{\sqrt{2}}(r+i)$. Using the values of $\beta$ from (\ref{2beta}) we obtain that $\frac{\beta}{\sqrt{2}}(r+i)=
-\frac{1}{2}(r-1+i(r+1))$, $\frac{1}{2}(-r-1+i(r-1))$, $\frac{1}{2}(r+1+i(-r+1))$, or $\frac{1}{2}(r-1+i(r+1))$.  So the values of $\Re e p_{31}\in{\bf Z}$ are $\frac{1-r}{2}$, $-\frac{r+1}{2}$, $\frac{r+1}{2}$, or $\frac{r-1}{2}$ and therefore $r$ is an odd integer.

We recall that 
\begin{eqnarray}
Q(G) & \geq &  
\delta^{2}\left(\frac{r}{\delta^{2}}+
\frac{1}{\delta^{2}}\Re e(p_{32}\bar{p}_{33})x_{1}- 
\frac{1}{\delta^{2}}\Im  m(p_{32}\bar{p}_{33})y_{1}+x_{2}\right)^{2}+ \nonumber \\
\ & \ & \delta^{2}\left(\frac{\frac{|\gamma|^{2}}{2}}{\delta^{2}}+
\frac{1}{\delta^{2}}\Im m(p_{32}\bar{p}_{33})x_{1}+ 
\frac{1}{\delta^{2}}\Re e(p_{32}\bar{p}_{33})y_{1}+y_{2}\right)^{2}. 
\label{e118}
\end{eqnarray}
We will only consider the case  $g_{12}\bar{g}_{13}=-2i$ because the other case can be handled similarly.   Since $p_{32}\bar{p}_{33}=-2i$ the right-hand side of (\ref{e118}) is equal to 
\begin{equation} 
2(\frac{r}{2}+y_{1}+x_{2})^{2}+2(\frac{1}{2}-x_{1}+y_{2})^{2}.
\label{e121}
\end{equation}
We may assume that $\frac{1+\sqrt{3}-\sqrt{2}}{2}\leq x_{1}\leq 1$. Indeed, otherwise we have $0\leq x_{1}\leq \frac{1+\sqrt{3}-\sqrt{2}}{2}$ and 
$$\frac{1}{2}-x_{1}+y_{2}\geq \frac{1}{2}-\frac{1+\sqrt{3}-\sqrt{2}}{2}+\frac{\sqrt{3}}{2}=\frac{1}{\sqrt{2}}.$$
According to (\ref{e121}) we get the lower bound $1$ for all $z\in{\cal F}_{1}$. 

We know that $0\leq y_{1}\leq 1-x_{1} \leq \frac{1-\sqrt{3}+\sqrt{2}}{2}$. Therefore 
$$\frac{|r|}{2}\leq |\frac{r}{2}+y_{1}+x_{2}|+y_{1}+|x_{2}|\leq |\frac{r}{2}+y_{1}+x_{2}|+\frac{1-\sqrt{3}+\sqrt{2}}{2}+\frac{1}{2}.$$
So if $|r|\geq 3$ then 
$$|\frac{r}{2}+y_{1}+x_{2}|\geq \frac{|r|}{2}-1+\frac{\sqrt{3}-\sqrt{2}}{2}\geq \frac{3}{2}-1+\frac{\sqrt{3}-\sqrt{2}}{2}=\frac{1+\sqrt{3}-\sqrt{2}}{2}>0.$$
Using the trivial lower bound $\frac{1}{2}-x_{1}+y_{2}\geq \frac{1}{2}-1+\frac{\sqrt{3}}{2}$ and the convexity of $x^{2}$ we get 
\begin{eqnarray*}
Q(G) & = & 2(\frac{r}{2}+y_{1}+x_{2})^{2}+2(\frac{1}{2}-x_{1}+y_{2})^{2}\\
 & \geq &  2(\frac{1+\sqrt{3}-\sqrt{2}}{2})^{2}+2(\frac{\sqrt{3}-1}{2})^{2}\\
\ & = & \frac{(1+\sqrt{3}-\sqrt{2})^{2}+(\sqrt{3}-1)^{2}}{2}\\
\ & \geq & (\frac{2\sqrt{3}-\sqrt{2}}{2})^{2}\geq 1
\end{eqnarray*}
for all $z\in {\cal F}_{1}$. This completes the proof of (ii). 
According to Lemma \ref{l121} and \ref{l124} the possible values of $r$ are odd integers between $-2$ and $2$. This completes the proof of Lemma \ref{l124}.  

\vskip 0.3cm
The remaining two cases are $g_{12}\bar{g}_{13}=-2i,-2$, with $r=-1,1$, $|\gamma|^{2}=2$, and $\beta$ as in (\ref{2beta}). 

Consider first the case $g_{12}\bar{g}_{13}=-2i$. In this case $g_{13}=ip_{33}=i\sqrt{2}\beta$, $g_{11}=ip_{31}=i\frac{\beta}{\sqrt{2}}(r+i)$, $g_{12}=ip_{32}=\sqrt{2}\beta$. So $Q(G)=|i\frac{\beta}{\sqrt{2}}(r+i)+\sqrt{2}\beta z_{1}+i\sqrt{2}\beta z_{2}|^{2}$, and therefore 

$$|\det G'(z)|^{2}=|\frac{r+i}{\sqrt{2}}-i\sqrt{2} z_{1}+
\sqrt{2}z_{2}|^{-6},\ \ r=-1,1.$$

The stability component of $G$ can be computed easily by using the Langlands
decomposition: 

$$P_{5+\frac{1+r}{2}}=\left(\begin{array}{ccc}
\frac{\beta}{\sqrt{2}} & 0 & 0 \\
-\beta^{-2} & \beta^{-2} & 0 \\ 
\frac{\beta}{\sqrt{2}}(r+i) &  -i\sqrt{2}\beta & \sqrt{2}\beta 
\end{array}\right)$$
with $r=-1,1$. Select $\beta=\frac{-1+ri}{\sqrt{2}}$. Consider the transformation 
$G_{5+\frac{1+r}{2}} \equiv N_{(1,\frac{r+i}{2})}\circ J\circ P_{5+\frac{1+r}{2}}$. This transformation is in the Picard modular group since its entries are Gaussian integers. The Jacobian determinant of  $G_{5+\frac{1+r}{2}}$ is 
$$|\det G_{5+\frac{1+r}{2}}'(z)|^{2}=|\frac{r+i}{\sqrt{2}}-i\sqrt{2} z_{1}+\sqrt{2}z_{2}|^{-6},\ \ r=-1,1.$$
The transformation 
$G_{5+\frac{1+r}{2}}$ satisfies the Siegel property, that is, there a point  
$R_{5+\frac{1+r}{2}}\in S(L)$ such that $G_{5+\frac{1+r}{2}}(R_{5+\frac{1+r}{2}})\in S(L)$. Indeed, since   $G_{5+\frac{1+r}{2}}(z_{1}, iz_{1})=(z_{1}, iz_{1})$ one can select $R_{5+\frac{1+r}{2}}=(1-2\epsilon+i\epsilon, -\epsilon+i(1-2\epsilon))$, $r=-1,1$, for some small $\epsilon>0$.

Now consider the second case $g_{12}\bar{g}_{13}=-2$. In this case $g_{13}=ip_{33}=i\sqrt{2}\beta$, $g_{11}=ip_{31}=i\frac{\beta}{\sqrt{2}}(r+i)$, $g_{12}=ip_{32}=-i\sqrt{2}\beta$. So $Q(G)=|i\frac{\beta}{\sqrt{2}}(r+i)-i\sqrt{2}\beta z_{1}+i\sqrt{2}\beta z_{2}|^{2}$, and therefore 

$$|\det G'(z)|^{2}=|\frac{r+i}{\sqrt{2}}-\sqrt{2} z_{1}+
\sqrt{2}z_{2}|^{-6},\ \ r=-1,1.$$

The stability component of $G$ can be computed easily by using the Langlands
decomposition: 

$$P_{7+\frac{1+r}{2}}=\left(\begin{array}{ccc}
\frac{\beta}{\sqrt{2}} & 0 & 0 \\
-i\beta^{-2} & \beta^{-2} & 0 \\ 
\frac{\beta}{\sqrt{2}}(r+i) &  -\sqrt{2}\beta & \sqrt{2}\beta 
\end{array}\right)$$
with $r=-1,1$. Select $\beta=\frac{-1+ri}{\sqrt{2}}$. Consider the transformation 
$G_{7+\frac{1+r}{2}} \equiv N_{(i,\frac{r+i}{2})}\circ J\circ P_{7+\frac{1+r}{2}}$. This transformation is in the Picard modular group since its entries are Gaussian integers. The Jacobian determinant of  $G_{7+\frac{1+r}{2}}$ is 
$$|\det G_{7+\frac{1+r}{2}}'(z)|^{2}=|\frac{r+i}{\sqrt{2}}-\sqrt{2} z_{1}+\sqrt{2}z_{2}|^{-6},\ \ r=-1,1.$$
The transformation 
$G_{7+\frac{1+r}{2}}$ satisfies the Siegel property, that is, there is a point  
$R_{7+\frac{1+r}{2}}\in S(L)$ such that $G_{7+\frac{1+r}{2}}(R_{7+\frac{1+r}{2}})\in S(L)$ because    $G_{7+\frac{1+r}{2}}(z_{1}, z_{1})=(z_{1}, z_{1})$.

\subsubsection{Class of automorphisms with dilation $\delta=2$}
Let $G=[g_{jk}]\in\Gamma\setminus\Gamma_{\infty}$ be a holomorphic automorphism of ${\bf CH}^{2}$ with dilation parameter  $\delta=|g_{13}|=2$. Then it follows from 
(\ref{e101}) and (\ref{e1015}) that $|p_{33}|^{2}=\delta^{2}=4$. This determines the possible values of the Gaussian integer $p_{33}=\delta\beta$ and $\beta$:
\begin{equation}
p_{33}=-2, \  2, -2i, \ 2i, \  \  \beta=-1,\ 1,\ -i,\ i.
\label{3beta}
\end{equation}

Since  $\beta^{2}=1, -1$ using (\ref{e101}) again we get from $\bar{\gamma}=-i\beta^{2}p_{32}$ that $\gamma$ is a Gaussian integer. Since $|\gamma|^{2}=2\Im m (p_{31}\bar{p}_{33})\in 2{\bf Z}$ and $0<|\gamma|^{2}<\gamma_{0}(2)^{2}<3.1^{2}<10$ according to Remark \ref{r103} and Lemma \ref{l112}, the possible values of $|\gamma|^{2}$ are $2$, $4$ and $8$. Indeed $6$ is not a possible value since it can not be written as sum of squares of two integers.   

Our next step is to eliminate the values $2$ and $8$ for $|\gamma|^{2}$.

\begin{lem}
Let $G=[g_{jk}]\in\Gamma\setminus\Gamma_{\infty}$ be a holomorphic automorphism of ${\bf CH}^{2}$ with $|g_{13}|=|ip_{33}|=\delta=2$ and $|g_{12}|^{2}=|\gamma|^{2}=2, 8$. Then  $|g_{11}+g_{12}z_{1}+g_{13}z_{2}|^{2}\geq 1$ for all $z\in {\cal F}_{1}$.
\label{l201}
\end{lem}

\noindent{\bf Proof of Lemma \ref{l201}.} Consider the case $|\gamma|^{2}=2$. We recall that 
\begin{eqnarray}
Q(G) & \geq &  
\delta^{2}\left(\frac{\frac{|\gamma|^{2}}{2}}{\delta^{2}}+
\frac{1}{\delta^{2}}\Im m(p_{32}\bar{p}_{33})x_{1}+ 
\frac{1}{\delta^{2}}\Re e(p_{32}\bar{p}_{33})y_{1}+y_{2}\right)^{2} \nonumber \\
\ & = & 4(\frac{1}{4}+\frac{1}{4}\Im m(p_{32}\bar{p}_{33})x_{1}+\frac{1}{4}\Re e(p_{32}\bar{p}_{33})y_{1}+y_{2})^{2}. 
\label{e201}
\end{eqnarray} 
 
Here $\Re e(p_{32}\bar{p}_{33})$, $\Im m(p_{32}\bar{p}_{33})$ are integers 
and  
$$|\Re e(p_{32}\bar{p}_{33})|^{2}+|\Im m(p_{32}\bar{p}_{33})|^{2}=|p_{32}\bar{p}_{33}|^{2}=
\delta^{2}|\gamma|^{2}=8.$$
 So the possible value of the pair 
$(|\Re e(p_{32}\bar{p}_{33})|,|\Im m(p_{32}\bar{p}_{33})|)$ 
is $(2,2)$. This gives the lower bound 
$\frac{1}{4}-\frac{1}{4}|\Im m(p_{32}\bar{p}_{33})|x_{1}-\frac{1}{4}|\Re e(p_{32}\bar{p}_{33})|y_{1}+y_{2}=\frac{1}{4}-\frac{1}{2}x_{1}-\frac{1}{2}y_{1}+y_{2}=\frac{1}{2}-\frac{1}{2}x_{1}-\frac{1}{2}y_{1}+y_{2}-\frac{1}{4}
\geq y_{2}-\frac{1}{4}\geq \frac{\sqrt{3}}{2}-\frac{1}{4}\geq \frac{1}{2}$ for the quantity in the parentheses on the right hand side of (\ref{e201}). 
Therefore 
$$Q(G)\equiv |g_{11}+g_{12}z_{1}+g_{13}z_{2}|^{2}\geq \delta^{2}
(\frac{1}{2})^{2}=1$$
 for all $z\in {\cal F}_{1}$. 

A similar argument works for the case $|\gamma|^{2}=8$. Indeed, in this case  $$|\Re e(p_{32}\bar{p}_{33})|^{2}+|\Im m(p_{32}\bar{p}_{33})|^{2}=|p_{32}\bar{p}_{33}|^{2}=
\delta^{2}|\gamma|^{2}=32,$$
 the possible value of the pair 
$(|\Re e(p_{32}\bar{p}_{33})|,|\Im m(p_{32}\bar{p}_{33})|)$ 
is $(4,4)$. Therefore
\begin{eqnarray*} 
Q(G) & \geq  & 4(1+\frac{1}{4}\Im m(p_{32}\bar{p}_{33})x_{1}+\frac{1}{4}\Re e(p_{32}\bar{p}_{33})y_{1}+y_{2})^{2} \\
 & \geq & 
4(1-x_{1}-y_{1}+y_{2})^{2}\geq 4y_{2}^{2}\geq 4(\frac{\sqrt{3}}{2})^{2}>1
\end{eqnarray*}
 for all $z\in {\cal F}_{1}$.

This completes the proof of Lemma \ref{l201}.

Now we turn our attention to the case $\delta=2$, $|\gamma|^{2}=4$. In this case 
$|\Re e(p_{32}\bar{p}_{33})|^{2}+|\Im m(p_{32}\bar{p}_{33})|^{2}=|p_{32}\bar{p}_{33}|^{2}=
\delta^{2}|\gamma|^{2}=16$. Then we know that the pair of integers 
$(|\Re e(p_{32}\bar{p}_{33})|, |\Im m(p_{32}\bar{p}_{33})|)$ is equal to $(0,4)$ or $(4,0)$ so $p_{32}\bar{p}_{33}=-4i, 4i, -4, 4$. 
The next lemma further reduces the possible values of $p_{32}\bar{p}_{33}$.

\begin{lem}
Let $G=[g_{jk}]\in\Gamma\setminus\Gamma_{\infty}$ be a holomorphic automorphism of ${\bf CH}^{2}$ with $|g_{13}|=|ip_{33}|=\delta=2$, $|g_{12}|^{2}=|\gamma|^{2}=4$ and $g_{12}\bar{g}_{13}=4i, 4$. Then  $|g_{11}+g_{12}z_{1}+g_{13}z_{2}|^{2}\geq 1$ for all $z\in {\cal F}_{1}$.
\label{l205}
\end{lem}

\noindent{\bf Proof of Lemma \ref{l205}.} Similarly to (\ref{e115}) we have 
$$|g_{11}+g_{12}z_{1}+g_{13}z_{2}|^{2} \geq   
4(\frac{1}{2}+\frac{1}{4}\Im m(p_{32}\bar{p}_{33})x_{1}+\frac{1}{4}\Re e(p_{32}\bar{p}_{33})y_{1}+y_{2})^{2}$$
in this case. If $g_{12}\bar{g}_{13}=p_{32}\bar{p}_{33}=4i$ then we get the lower bound 
$$4(\frac{1}{2}+x_{1}+y_{2})^{2}\geq 4(\frac{1}{2})^{2}= 1.$$
The other case, $g_{12}\bar{g}_{13}=p_{32}\bar{p}_{33}=4$,  gives the lower bound 
$$4(\frac{1}{2}+y_{1}+y_{2})^{2}\geq 4(\frac{1}{2})^{2}= 1.$$
These lower bounds imply Lemma \ref{l205}.

\begin{lem}
Let $G=[g_{jk}]\in\Gamma\setminus\Gamma_{\infty}$ be a holomorphic automorphism of ${\bf CH}^{2}$ with $|g_{13}|=|ip_{33}|=\delta=2$, $|g_{12}|^{2}=|\gamma|^{2}=4$ and $g_{12}\bar{g}_{13}=-4i, -4$. 

(i) If $|r|\equiv|\Re e g_{11}\bar{g}_{13}|\geq 4$ or $r=0$ then  $Q(G)\geq 1$ for all $z\in {\cal F}_{1}$.

(ii) The number $r=\Re e g_{11}\bar{g}_{13}\in{\bf Z}$ is even. 

(iii) The possible values of $r\in{\bf Z}$ are $\pm 2$. 

\label{l210}
\end{lem}

\noindent{\bf Proof of Lemma \ref{l210}}. To prove (i) we use (\ref{e103}) and 
(\ref{e109}). It is easy to see that 
\begin{equation}
Q(G)=\left\{\begin{array}{ll}
(\frac{r}{2}+2y_{1}+2x_{2})^{2}+(1-2x_{1}+2y_{2})^{2} & \mbox{ if } g_{12}\bar{g}_{13}=-4i, \\
(\frac{r}{2}-2x_{1}+2x_{2})^{2}+(1-2y_{1}+2y_{2})^{2} & \mbox{ if } g_{12}\bar{g}_{13}=-4 .
\end{array}
\right. 
\label{e205}
\end{equation}
We will concentrate on the first case only, the other case is similar. We may assume that $\sqrt{3}/2\leq x_{1}\leq 1$. Indeed if $0\leq x_{1}\leq \sqrt{3}/2$ then $1-2x_{1}+2y_{2}\geq 1-2x_{1}+\sqrt{3}\geq 1$ so $Q(G)\geq 1$. We know that $0\leq y_{1}\leq 1-x_{1}\leq 1 -\sqrt{3}/2$ in ${\cal F}_{1}$. This means $|r|/2\leq |\frac{r}{2}+2y_{1}+2x_{2}|+2y_{1}+2|x_{2}|\leq 
|\frac{r}{2}+2y_{1}+2x_{2}|+2-\sqrt{3}+1$. If $|r|\geq 4$ then we have 
$$0<2+\sqrt{3}-3\leq |\frac{r}{2}|+\sqrt{3}-3\leq |\frac{r}{2}+2y_{1}+2x_{2}|.$$
Therefore $Q(G)\geq (\sqrt{3}-1)^{2}+(-1+\sqrt{3})^{2}>1$ in ${\cal F}_{1}.$ 

 If $r=0$ then 
$$Q(G)=4[(y_{1}+x_{2})^{2}+(\frac{1}{2}-x_{1}+y_{2})^{2}].$$
We again may assume that $\sqrt{3}/2\leq x_{1}\leq 1$. Since $|x_{2}|\leq \frac{1}{2}$, $0\leq x_{1}\leq 1$, we have $\sqrt{1-x_{2}^{2}}\geq \sqrt{3}/2$ and $\frac{1}{2}-x_{1}+y_{2}\geq 
\frac{1}{2}-x_{1}+\sqrt{1-x_{2}^{2}}\geq \frac{1}{2}-x_{1}+\sqrt{3}/2\geq (\sqrt{3}-1)/2>0$ in ${\cal F}_{1}$. Therefore 
\begin{eqnarray*}
Q(G) & \geq & 4[(y_{1}+x_{2})^{2}+(\frac{1}{2}-x_{1}+\sqrt{1-x_{2}^{2}})^{2}] \\
 & = & 4[y_{1}^{2}+2y_{1}x_{2}+x_{2}^{2}+\frac{1}{4}-x_{1}+x_{1}^{2}+2(\frac{1}{2}-x_{1})\sqrt{1-x_{2}^{2}}+1-x_{2}^{2}] \\
 & \geq & 4[\frac{1}{4}-x_{1}-2|x_{2}|y_{1}+1+x_{1}^{2}+2(\frac{1}{2}-x_{1})\sqrt{1-x_{2}^{2}}].
\end{eqnarray*}
Since $1/2<\sqrt{3}/{2}\leq x_{1}$ so $(\frac{1}{2}-x_{1})\sqrt{1-x_{2}^{2}}\geq \frac{1}{2}-x_{1}$. Therefore

\begin{eqnarray*}
Q(G) & \geq & 4[\frac{1}{4}-x_{1}-y_{1}+1+x_{1}^{2}+2(\frac{1}{2}-x_{1})]\\
 & \geq & 4[\frac{1}{4}+(x_{1}-1)^{2}]\geq 1  
\end{eqnarray*}
in ${\cal F}_{1}$. 

(ii) It follows from (\ref{e101}), (\ref{e1015}) and (\ref{3beta}) that $g_{11}=\beta(\frac{r}{2}i-1)$ with $\beta=\pm 1,\pm i$. Since $g_{11}$ is a Gaussian integer, we get that $\frac{r}{2}\in{\bf Z}$, that is $r$ is even. 

(iii) Since $0<|r|<4$, $r\in{\bf Z}$ is even, the only possible value of $r$ is $\pm 2$. 
This completes the proof of Lemma \ref{l210}.

Now we summarize the case $\delta=2$. Let $g\in\Gamma\setminus\Gamma_{\infty}$ be an automorphism with $\delta=|g_{13}|=2$. Then it follows from Lemma  \ref{l201}, \ref{l205} and \ref{l210} that either $|\det G'(z)|\leq 1$ for all $z\in{\cal F}_{1}$, or $g_{11}=i\beta(\frac{r}{2}+i)$, $g_{12}=2\beta, -2i\beta$, $g_{13}=2i\beta$ with $\beta=\pm 1,\pm i$ and $r=\pm 2$. In the former case the transformation $G$ does not contribute to the fundamental domain ${\cal F}$. In the latter case 

\begin{equation*}
Q(G)=\left\{\begin{array}{ll}
|i\beta(\frac{r}{2}+i)+2\beta z_{1}+2i\beta z_{2}|^{2}=|\frac{r}{2}+i-2i z_{1}+2z_{2}|^{2}
 & \mbox{ if } g_{12}=2\beta, \\
|i\beta(\frac{r}{2}+i)-2i\beta z_{1}+2i\beta z_{2}|^{2}=|\frac{r}{2}+i-2 z_{1}+2z_{2}|^{2}
 & \mbox{ if } g_{12}=-2i\beta
\end{array}
\right. 
\end{equation*}
with $r=\pm 2$. We recall that $Q(G_{5+\frac{1+k}{2}})=\frac{1}{2}|k+i-2iz_{1}+2z_{2}|^{2}$,  
$Q(G_{7+\frac{1+k}{2}})=\frac{1}{2}|k+i-2z_{1}+2z_{2}|^{2}$ with $k=\pm 1$. So $Q(G)$ is a multiple of one of $Q(G_{j})$ for $j=5,6,7,8$. Indeed 
\begin{equation*}
Q(G)=\left\{\begin{array}{ll}
2Q(G_{5}) & \mbox{ if } g_{12}=2\beta, \ r=-2 \\
2Q(G_{6}) & \mbox{ if } g_{12}=2\beta, \ r=2 \\
2Q(G_{7}) & \mbox{ if } g_{12}=-2i\beta, \ r=-2 \\
2Q(G_{8}) & \mbox{ if } g_{12}=-2i\beta, \ r=2. 
\end{array}
\right. 
\end{equation*}
If for some point $z\in{\cal F}_{1}$, $Q(G)(z)<1$, then $Q(G_{j})<\frac{1}{2}$ for a $j\in\{5,6,7,8\}$. 
Therefore no transformation with dilation parameter $\delta=2$ contributes to the fundamental domain ${\cal F}$.

\noindent{\bf Proof of Theorem \ref{maint05}}. The proof follows the outline described in section 2.1. Let $G$ be one of the transformations $G_{2},\dots,G_{N}$ appearing in the description of $\cal F$ in (\ref{sef}). According to Lemma \ref{l103} and \ref{l115} the dilation parameter of $G$ is $1$, $\sqrt{2}$, or $2$. In section 3.2.1  we proved that if the dilation parameter of $G$ is $1$ then either $G\in\{G_{2}, G_{3}, G_{4}\}$ or $|\det G'(z)|\leq 1$ for all $z\in{\cal F}_{1}\supset{\cal F}$. In the latter case the presence of the condition $|\det G'(z)|\leq 1$ in (\ref{sef}) is not relevant. Similarly if the dilation parameter of $G$ is $\sqrt{2}$, we proved in section 3.2.2 that either $G\in\{G_{5}, G_{6}, G_{7}, G_{8}\}$ or $|\det G'(z)|\leq 1$ for all $z\in{\cal F}_{1}\supset{\cal F}$. In the latter case the transformation $G$ does not contribute to the description of $\cal F$ in (\ref{sef}). If the dilation parameter of $G$ is $\delta =2$, we proved in section 3.2.3 that either $|\det G'(z)|\leq 1$ for all $z\in{\cal F}_{1}\supset{\cal F}$ again, or there is a  transformation $\tilde{G}\in\{G_{5}, G_{6}, G_{7}, G_{8}\}$ such that $|\det G'(z)|\leq |\det \tilde{G}'(z)|$ for all $z\in{\cal F}_{1}\supset{\cal F}$. Therefore the restriction $|\det G'(z)|\leq 1$ in (\ref{sef}) does not have an effect on the fundamental domain for any transformation $G\in\{G_{2},\dots,G_{N}\}$ with dilation parameter $\delta=2$.  
 This completes the proof of Theorem \ref{maint05}.

\noindent{\bf Proof of Theorem \ref{maint1}}. Theorem \ref{maint1} follows from Theorem \ref{maint05}. Since the transformations $G_{2},\dots, G_{8}$ are explicitly given, one can easily compute the Jacobian determinants $\det G'_{j}(z)$, $j=2,\dots, 8$ according to (\ref{ejacobi}). The inequalities $|\det G'_{j}(z)|\leq 1$ $j=2,\dots ,8$ in the statement of Theorem \ref{maint1} are identical to the last three inequalities appearing in the definition of ${\cal F}$ in Theorem \ref{maint05}.  

All eight transformations $G_{1},\dots, G_{8}$ are needed in the description of the fundamental domain in ${\cal F}$ in Theorems \ref{maint05}, \ref{maint1}. In other words, if one of $G_{1},\dots, G_{8}$ is omitted in the description of $\cal F$, the set obtained is strictly larger than $\cal F$. This can be easily verified by examining the location of the points 
$$P_{1}(0.2+0.2i, 0.1+0.8i),\ P_{2}(0.25+0.72i, 0.45+0.9i),$$ $$P_{3}(0.39+0.6i, -0.45+0.9i),\ P_{4}(0.72+0.25i, -0.45+0.9i),$$
$$P_{5}(0.7+0.05i, 0.45+0.9i),\ P_{6}(0.7+0.02i, -0.45+0.9i),$$
$$P_{7}(0.02+0.7i, 0.45+0.9i),\ P_{8}(0.02+0.7i, -0.45+0.9i).$$
Simple calculation shows, that, if $1\leq k\leq 4$ then 
\begin{equation*}
Q_{j}(P_{k})=\left\{\begin{array}{ll}  
\geq 1 & \mbox{ if } 1\leq j\leq 4,\ j\not=k,\\
<1 & \mbox{ if } j=k,\\
\geq 2 & \mbox{ if } 5\leq j\leq 8;
\end{array}
\right.
\end{equation*}
and if $5\leq k\leq 8$ then 
\begin{equation*}
Q_{j}(P_{k})=\left\{\begin{array}{ll}  
\geq 1 & \mbox{ if } 1\leq j\leq 4,\\
<2 & \mbox{ if } j=k,\\
\geq 2 & \mbox{ if } 5\leq j\leq 8,\ j\not=k.
\end{array}
\right.
\end{equation*}
Moreover the points $P_{j}$, $j=1,\dots, 8$, satisfy the conditions $0\leq\Re e z_{1}$, $0\leq \Im m z_{1}$, $\Re e z_{1}+\Im m z_{1}\leq 1$, $|\Re e z_{2}|\leq\frac{1}{2}$.
The proof of Theorem \ref{maint1} has been finished.   

\noindent{\bf Proof of Theorem \ref{maint3}}. (i) The transformation $S$ preserves the four inequalities,  $0\leq\Re e z_{1}$, $0\leq\Im m z_{1}$, $\Re e z_{1}+\Im m z_{1}\leq 1$, and $|\Re e z_{2}|\leq\frac{1}{2}$ in the definition of $\cal F$ automatically. An easy calculation shows that 
$$Q_{1}(S(z))=Q_{1}(z),$$  
$$Q_{2}(S(z))=Q_{4}(z),$$ 
$$Q_{3}(S(z))=Q_{3}(z),$$ 
$$Q_{4}(S(z))=Q_{2}(z),$$
$$Q_{5}(S(z))=Q_{8}(z),$$  
$$Q_{6}(S(z))=Q_{7}(z),$$ 
$$Q_{7}(S(z))=Q_{6}(z),$$
$$Q_{8}(S(z))=Q_{5}(z).$$ 
Therefore the transformation $S$ maps $\cal F$ onto itself since $S^{2}=I$.

(ii) It is enough to show that the inequalities $|\Re e z_{2}|\leq \frac{1}{2}$, $\Im m z_{2}>0$, $|z_{2}|\geq 1$ imply that $(0,z_{2})\in{\cal F}_{1}$, that is,  
$$|r+i+z_{2}|^{2}\geq 1,\ \mbox{ for }\ r=-1,0,1$$
and 
$$|r+i+2z_{2}|^{2}\geq 2,\ \mbox{ for }\ r=-1,1.$$
This follows from $\Im m z_{2}\geq\sqrt{3}/2$, $|r+i+z_{2}|^{2}\geq (1+\Im m z_{2})^{2}\geq 1$ and  
$|r+i+2z_{2}|^{2}\geq (1+2\Im m z_{2})^{2}\geq (1+\sqrt{3})^{2}\geq 2$. 

(iii) We will prove the statement for any $a\geq 1+2\sqrt{2}$. If $\Im m z_{2}\geq a$, $|z_{1}|\leq 1$, $|r|\leq 1$ then 
\begin{eqnarray*}
1+2\sqrt{2}\leq a\leq |z_{2}| & \leq & |r+i-(1+i)z_{1}+z_{2}|+|r+i-(1+i)z_{1}|\\
\ & \leq & |r+i-(1+i)z_{1}+z_{2}|+|r+i|+|1+i|\\
\ & \leq & |r+i-(1+i)z_{1}+z_{2}|+2\sqrt{2}.
\end{eqnarray*}
So $|r+i-(1+i)z_{1}+z_{2}|\geq 1$. 
Similarly $|r+i-2iz_{1}+2z_{2}|\geq \sqrt{2}$ and $|r+i-2z_{1}+2z_{2}|\geq \sqrt{2}$. Therefore $\{z\in{\bf C}^{2};\ \ 0\leq\Re e z_{1},\ 0\leq\Im m z_{1},\ \Re e z_{1}+\Im m z_{1}\leq 1,\ |\Re e z_{2}|\leq 1/2,\ \Im m z_{2}\geq a\}\subset{\cal F}$. This completes the proof of Theorem \ref{maint3}.

\end{document}